# The two-parameter Poisson–Dirichlet point process

## KENJI HANDA

*Department of Mathematics, Saga University, Saga 840-8502, Japan.*
*E-mail: handa@ms.saga-u.ac.jp*

The two-parameter Poisson–Dirichlet distribution is a probability distribution on the totality of positive decreasing sequences with sum 1 and hence considered to govern masses of a random discrete distribution. A characterization of the associated point process (that is, the random point process obtained by regarding the masses as points in the positive real line) is given in terms of the correlation functions. Using this, we apply the theory of point processes to reveal the mathematical structure of the two-parameter Poisson–Dirichlet distribution. Also, developing the Laplace transform approach due to Pitman and Yor, we are able to extend several results previously known for the one-parameter case. The Markov–Krein identity for the generalized Dirichlet process is discussed from the point of view of functional analysis based on the two-parameter Poisson–Dirichlet distribution.

## 1. Introduction

This paper is concerned with a two-parameter family of probability distributions on the infinite-dimensional simplex $\nabla_\infty$ of non-negative decreasing sequences with sum 1,

$$\nabla_\infty = \left\{ (v_i) = (v_1, v_2, \ldots) : v_1 \geq v_2 \geq \cdots \geq 0, \sum v_i = 1 \right\}.$$

It extends the one-parameter family of distributions known as *Poisson–Dirichlet distributions*, which were introduced by Kingman [32]; see, for example, [1, 44, 50] for related topics and bibliographic information. Pitman and Yor [45] defined the two-parameter Poisson–Dirichlet distribution, denoted $\mathrm{PD}(\alpha, \theta)$, in the following manner. Given $0 \leq \alpha < 1$ and $\theta > -\alpha$, define a sequence $(\widetilde{V}_i)$ of random variables by

$$\widetilde{V}_1 = Y_1, \qquad \widetilde{V}_i = (1 - Y_1) \cdots (1 - Y_{i-1}) Y_i \qquad (i = 2, 3, \ldots), \tag{1.1}$$

where $Y_1, Y_2, \ldots$ are independent and each $Y_i$ is Beta$(1 - \alpha, \theta + i\alpha)$-distributed. Let $(V_i)$ be the decreasing order statistics of $(\widetilde{V_i})$, namely $V_1 \geq V_2 \geq \cdots$ are the ranked values of $(\widetilde{V_i})$, and define PD$(\alpha, \theta)$ to be the law of $(V_i)$ on $\nabla_\infty$. There is some background to the study of these distributions, as explained in [45] and [44]. In particular, a number of results concerning PD$(\alpha, \theta)$ were obtained in [45]. The proofs there use auxiliary random variables and related processes (such as stable subordinators and gamma processes) and require deep insight into them. As for the original Poisson–Dirichlet distributions, which form a one-parameter family $\{\text{PD}(0, \theta) : \theta > 0\}$ and correspond to the gamma processes, certain independence property often makes the analysis relatively easier. See [50] and [51] for extensive discussions.

One purpose of this article is to provide another approach to study the two-parameter family $\{\text{PD}(\alpha, \theta) : 0 \leq \alpha < 1, \theta > -\alpha\}$, based mainly on conventional arguments in the theory of point processes; see, for example, [7] for general accounts of the theory. This means that a random element $(V_i)$ governed by PD$(\alpha, \theta)$ is studied through the random point process $\xi := \sum \delta_{V_i}$, which we call the (two-parameter) Poisson–Dirichlet point process with parameters $(\alpha, \theta)$, or simply the PD$(\alpha, \theta)$ process. Note that although the above $\xi$ is a 'non-labeled' object, it is sufficient to recover the law of the ranked sequence $(V_i)$. For example, we have $\xi([t, \infty)) = 0$ precisely when $V_1 < t$. More generally, for each $n = 2, 3, \ldots$, quantitative information on $V_1, \ldots, V_n$ can be derived from $\xi$ by the principle of inclusion-exclusion. Among the many ways of characterizing a point process, we choose the one prescribed in terms of correlation functions. For each positive integer $n$, the $n$th correlation function of a point process is informally defined as the mean density of tuples of $n$ distinct points in the process. As far as the one-parameter family $\{\text{PD}(0, \theta) : \theta > 0\}$ is concerned, the idea of such an approach is not new. The correlation functions of the PD$(0, \theta)$ process were computed by Watterson [54], where these are referred to as the *multivariate frequency spectra*. (See also [39].)

One advantage of this approach can be described as follows. As Kingman [33] mentioned, the Poisson–Dirichlet distribution is "rather less than user-friendly". In other words, PD$(0, \theta)$ is not so easy to handle directly. In contrast, the correlation functions of the PD$(0, \theta)$ process obtained in [54] are of a certain simple form. By using them, Griffiths [19] obtained several distributional results for PD$(0, \theta)$. As a two-parameter generalization, we will find ((2.3) below) the correlation functions of the PD$(\alpha, \theta)$ process. This result will play a key role in revealing the mathematical structure of the two-parameter Poisson–Dirichlet distribution via the associated point process. Indeed, by exploiting some general tools from the theory of point processes, it will be possible to deduce results which extend previously known results for the one-parameter case, for example, the joint probability density obtained by Watterson [54] and the moment formula due to Griffiths [18]. It should be noted that Pitman and Yor [45] essentially found the joint probability density in the two-parameter setting. However, we emphasize that our expression (5.13) in Theorem 5.4 is quite consistent with Watterson's formula and useful for the purpose of further analysis.

Another aspect of results presented in this article is a development of the Laplace transform method discussed in the last section of [45]. Roughly speaking, the main result ((3.8) in Theorem 3.2) in this direction states that the (suitably modified) Laplace transform of a 'probability generating functional' of the $PD(\alpha, \theta)$ process is connected with that of a certain Poisson point process on $(0, \infty)$ through some nonlinear function. This will provide a powerful tool, especially for the study of certain asymptotic problems regarding $PD(\alpha, \theta)$ as $\theta \to \infty$. Such problems have been studied in the case where $\alpha = 0$, motivated by the study of population genetics; see [18] for background information. We will show that the presence of a fixed parameter $\alpha \in (0, 1)$ does not affect the validity of results of such type, generalizing a limit theorem of Griffiths [18] and the central limit theorem obtained by Joyce, Krone and Kurtz [28]. In the proofs, Theorem 3.2 combined with the Laplace method is seen to be very effective in demonstrating such results.

The aforementioned structure of $PD(\alpha, \theta)$ will also be used to explain certain properties of the generalized Dirichlet process $\eta := \sum V_i \delta_{X_i}$, where $X_i$ are i.i.d. random variables independent of $(V_i)$. The Dirichlet process was originally introduced by Ferguson [17] as a prior distribution for nonparametric problems and it corresponds to the case where $(V_i)$ is governed by $PD(0, \theta)$. In such a context, many attempts have been made to find the exact distribution of the (random) mean $M := \sum X_i V_i$ of $\eta$. A key ingredient for this is the so-called Markov–Krein identity, which relates the law of $M$ to the common law $\nu$ of $X_i$. Its two-parameter generalization was obtained by Tsilevich [49]; see also [31, 50, 53]. Our study of the generalized Dirichlet process will be based on functional analysis of $PD(\alpha, \theta)$: we regard this process as a 'functional' of $PD(\alpha, \theta)$ with underlying parameter $\nu$. It will be seen that Theorem 3.2 explains the mathematical structure behind the Markov–Krein identity, a complementary result of which is presented. We would expect that such a point of view would help us to discuss further the generalized Dirichlet process.

The paper is organized as follows. In the next section, we describe the formula for the correlation function of the $PD(\alpha, \theta)$ process. Section 3 deals with the Laplace transform method, which will turn out to be another basic tool in the subsequent argument. In Section 4, we discuss the Dickman function and its two-parameter version which describes the law of $V_1$ with respect to $PD(\alpha, \theta)$. The main result here ((4.9) in Theorem 4.2) is an integral equation satisfied by the generalized Dickman function. Section 5 shows how the correlation functions combined with techniques from the theory of point processes yield distributional information on $PD(\alpha, \theta)$. In Section 6, we illustrate the utility of the Laplace transform method established in Section 3 by proving some asymptotic results on $PD(\alpha, \theta)$ with $\theta$ large. Section 7 is concerned with the generalized Dirichlet process, some of whose properties will reduce to the results obtained in previous sections.

## 2. The correlation function of the $PD(\alpha, \theta)$ process

The main subject of this section is the correlation function. We begin with some definitions and concepts from the theory of point processes [7]. (See also Section 2 in [27] for a comprehensive exposition of this material.) Let $\xi$ be a point process taking values in $\mathbf{R}$. For simplicity, suppose that $\xi$ is expressed as $\xi = \sum \delta_{X_i}$ for some random variables $X_1, X_2, \ldots$. Throughout what follows, we assume that $\xi$ is simple, in the sense that $X_i \neq X_j (i \neq j)$ a.s. For any positive integer $n$, the $n$th correlation measure (also called the $n$th factorial moment measure) of $\xi$, if it exists, is defined to be a $\sigma$-finite measure

$\mu_n$ such that for non-negative measurable functions $f$ on $\mathbf{R}^n$,

$$E\left[\sum_{i_1,\ldots,i_n(\neq)} f(X_{i_1},\ldots,X_{i_n})\right] = \int_{\mathbf{R}^n} f(x_1,\ldots,x_n)\mu_n(\mathrm{d}x_1\cdots\mathrm{d}x_n), \qquad (2.1)$$

where the sum is taken over $n$-tuples of distinct indices. In particular, $\mu_1$ is the mean measure of $\xi$ and it follows that $\mu_n(\mathrm{d}x_1\cdots\mathrm{d}x_n)$ is necessarily symmetric in $x_1,\ldots,x_n$. If $\mu_n$ has a density with respect to the $n$-dimensional Lebesgue measure, the density is called the $n$th *correlation function* of $\xi$.

In what follows, $\alpha$ and $\theta$ are such that $0 \leq \alpha < 1$ and $\theta > -\alpha$, unless otherwise mentioned. Notice that the PD$(\alpha,\theta)$ process $\xi = \sum \delta_{\widetilde{V}_i} = \sum \delta_{V_i}$, with $(\widetilde{V}_i)$ and $(V_i)$ defined as in the [Introduction](#), is simple. Denote by $c_{n,\alpha,\theta}$ the product

$$\prod_{i=1}^n \frac{\Gamma(\theta+1+(i-1)\alpha)}{\Gamma(1-\alpha)\Gamma(\theta+i\alpha)} = \begin{cases} \theta^n, & \alpha=0, \\ \dfrac{\Gamma(\theta+1)\Gamma(\theta/\alpha+n)\alpha^{n-1}}{\Gamma(\theta+\alpha n)\Gamma(\theta/\alpha+1)\Gamma(1-\alpha)^n}, & 0<\alpha<1. \end{cases}$$

In view of the left-hand side, it is clear that

$$c_{m+n,\alpha,\theta} = c_{m,\alpha,\theta}c_{n,\alpha,\theta+\alpha m} \qquad (m,n\in\{0,1,2,\ldots\}) \qquad (2.2)$$

with the convention that $c_{0,\alpha,\theta} = 1$. We also need the following notation for the $n$-dimensional unit simplex:

$$\Delta_n = \{(v_1,\ldots,v_n)\colon v_1\geq 0,\ldots,v_n\geq 0, v_1+\cdots+v_n\leq 1\}.$$

In general, the indicator function of a set (or an event) $E$ is denoted by $\mathbf{1}_E$. The main result of this section generalizes Watterson's formula [54] ((18), page 644) for the correlation functions of the PD$(0,\theta)$ process to the PD$(\alpha,\theta)$ process.

**Theorem 2.1.** *For each $n = 1,2,\ldots,$ the $n$th correlation function of the PD$(\alpha,\theta)$ process is given by*

$$q_{n,\alpha,\theta}(v_1,\ldots,v_n) := c_{n,\alpha,\theta}\prod_{i=1}^n v_i^{-(\alpha+1)}\left(1-\sum_{j=1}^n v_j\right)^{\theta+\alpha n-1}\mathbf{1}_{\Delta_n}(v_1,\ldots,v_n). \qquad (2.3)$$

We will give a proof based on a known property of the size-biased permutation for PD$(\alpha,\theta)$. For the one-parameter family $\{\text{PD}(0,\theta)\colon\theta>0\}$, such an idea to derive (2.3) is suggested in [39], Corollary 7.4. From the definition of size-biased permutation (see [45] or Section 4 of [3]), we can make the following observation. If a sequence $(V_i)$ of random variables such that $V_i > 0$ $(i=1,2,\ldots)$ and $\sum V_i = 1$ is given, then, for each $n = 1,2,\ldots,$ the $n$th correlation measure $\mu_n$ of $\sum \delta_{V_i}$ exists and

$$\mu_n(\mathrm{d}v_1\cdots\mathrm{d}v_n) = \prod_{i=1}^n\left\{v_i^{-1}\left(1-\sum_{j=1}^{i-1}v_j\right)\right\}\mathbf{1}_{\Delta_n}(v_1,\ldots,v_n)\mu_n^\sharp(\mathrm{d}v_1\cdots\mathrm{d}v_n), \qquad (2.4)$$

where $\mu_n^\sharp$ is the joint law of $(V_1^\sharp, \ldots, V_n^\sharp)$, the first $n$ components of the size-biased permutation $(V_i^\sharp)$ of $(V_i)$. In (2.4) and throughout, we adopt the convention that $\sum_{j=1}^0 \cdots = 0$.

**Proof of Theorem 2.1.** Let $(V_i)$ have $\mathrm{PD}(\alpha, \theta)$ distribution. A remarkable result obtained in [37, 41, 43] then tells us that the law of $(V_i^\sharp)$ coincides with the law of $(\widetilde{V}_i)$ in (1.1). Hence for any non-negative measurable function $f$ on $\mathbf{R}^n$, the expectation $E[f(V_1^\sharp, \ldots, V_n^\sharp)]$ is given by

$$\int_{[0,1]^n} \mathrm{d}y_1 \cdots \mathrm{d}y_n \, f(v_1, \ldots, v_n) c_{n,\alpha,\theta} \prod_{i=1}^n \{y_i^{-\alpha}(1 - y_i)^{\theta + i\alpha - 1}\}, \qquad (2.5)$$

where $v_i = v_i(y_1, \ldots, y_n) = y_i \prod_{j=1}^{i-1}(1 - y_j)$, with the convention that $\prod_{j=1}^0 \cdots = 1$. Note that the mapping $(y_1, \ldots, y_n) =: \mathbf{y} \mapsto (v_1, \ldots, v_n) =: \mathbf{v}$ from $[0,1]^n$ to $\Delta_n$ has the inverse

$$y_i = y_i(\mathbf{v}) = v_i \left(1 - \sum_{j=1}^{i-1} v_j\right)^{-1} \qquad (i = 1, \ldots, n). \qquad (2.6)$$

Therefore, (2.5) becomes

$$c_{n,\alpha,\theta} \int_{\Delta_n} \mathrm{d}v_1 \cdots \mathrm{d}v_n \left|\frac{\partial \mathbf{y}}{\partial \mathbf{v}}\right| f(\mathbf{v}) \prod_{i=1}^n \{y_i^{-\alpha}(1 - y_i)^{\theta + i\alpha - 1}\}. \qquad (2.7)$$

Observing from (2.6) that

$$\prod_{i=1}^n \{y_i^{-\alpha}(1 - y_i)^{\theta + i\alpha - 1}\} = \prod_{i=1}^n v_i^{-\alpha} \left(1 - \sum_{j=1}^n v_j\right)^{\theta + \alpha n - 1}, \qquad (2.8)$$

we have, by (2.5) and (2.7),

$$E[f(V_1^\sharp, \ldots, V_n^\sharp)] = \int_{\mathbf{R}^n} \mathrm{d}v_1 \cdots \mathrm{d}v_n \, f(\mathbf{v}) q_{n,\alpha,\theta}(\mathbf{v}) \prod_{i=1}^n \left\{v_i \left(1 - \sum_{j=1}^{i-1} v_j\right)^{-1}\right\}, \qquad (2.9)$$

where $q_{n,\alpha,\theta}$ is defined by (2.3). With the help of (2.4), this proves Theorem 2.1. $\qquad \square$

It is known that correlation functions appear in the expansion of the 'probability generating function' of a random point process $\sum \delta_{X_i}$; see Section 5 in [7] for general accounts. This functional is defined to be the expectation of an infinite product of the form $\prod g(X_i)$. For the sake of clarity, we provide the following definition of infinite products. Given a sequence $\{t_i\}$ of complex numbers, define

$$\prod_{i=1}^{\infty}(1+t_i) = \begin{cases} \displaystyle\lim_{n\to\infty}\prod_{i=1}^{n}(1+t_i) \in [1,\infty], & \text{if } t_i \geq 0 \ (i=1,2,\ldots), \\ \displaystyle\lim_{n\to\infty}\prod_{i=1}^{n}(1+t_i) \in \mathbf{C}, & \text{if the limit exists.} \end{cases}$$

In general, it holds that

$$\prod_{i=1}^{\infty}(1+|t_i|) = 1 + \sum_{n=1}^{\infty}\sum_{i_1<\cdots<i_n}|t_{i_1}|\cdots|t_{i_n}| = 1 + \sum_{n=1}^{\infty}\frac{1}{n!}\sum_{i_1,\ldots,i_n(\neq)}|t_{i_1}|\cdots|t_{i_n}|. \quad (2.10)$$

Also, $\prod_{i=1}^{\infty}(1+t_i)$ admits the corresponding expansion whenever $\prod_{i=1}^{\infty}(1+|t_i|)$ is finite. In what follows, $E_{\alpha,\theta}[\cdots]$ denotes the expectation with respect to $\mathrm{PD}(\alpha,\theta)$.

**Corollary 2.2.** *For any measurable function* $\phi\colon(0,1]\to\mathbf{C}$,

$$E_{\alpha,\theta}\left[\prod_{i=1}^{\infty}(1+|\phi(V_i)|)\right]$$
$$= 1 + \sum_{n=1}^{\infty}\frac{c_{n,\alpha,\theta}}{n!}\int_{\Delta_n}\prod_{i=1}^{n}\frac{|\phi(v_i)|}{v_i^{\alpha+1}}\left(1-\sum_{j=1}^{n}v_j\right)^{\theta+\alpha n-1}\mathrm{d}v_1\cdots\mathrm{d}v_n. \quad (2.11)$$

*If the above series converges, then*

$$E_{\alpha,\theta}\left[\prod_{i=1}^{\infty}(1+\phi(V_i))\right]$$
$$= 1 + \sum_{n=1}^{\infty}\frac{c_{n,\alpha,\theta}}{n!}\int_{\Delta_n}\prod_{i=1}^{n}\frac{\phi(v_i)}{v_i^{\alpha+1}}\left(1-\sum_{j=1}^{n}v_j\right)^{\theta+\alpha n-1}\mathrm{d}v_1\cdots\mathrm{d}v_n. \quad (2.12)$$

**Proof.** (2.11) is immediate from (2.10), (2.1) and (2.3) together. Convergence of the series in (2.11) allows us to show (2.12) by the dominated convergence theorem. $\square$

For instance, if $\phi\colon(0,1]\to\mathbf{C}$ is a measurable function with support contained in $[\varepsilon,1]$ for some $0<\varepsilon<1$, the series in (2.11) is easily verified to converge. Roughly speaking, the assertion of Corollary 2.2 is equivalent to that of Theorem 2.1 (cf. Propositions 2.2 and 2.3 in [27]).

# 3. The use of Laplace transforms

This section is intended to provide a general tool (Lemma 3.1) and to exploit it in the study of PD($\alpha, \theta$) processes. It contains a certain inversion formula for Laplace transforms. Interestingly, in spite of the generality of the formulation, the formula involves Dirichlet measures. This seems to show an aspect of the analytic importance of such measures.

**Lemma 3.1.** *Let $\beta \geq 0$ and $\delta \geq -\beta$. Let $R$ be a function defined at least on a neighborhood of $0$ in which $R$ is expressed as an absolutely convergent series of the form $R(u) = r_1 u + r_2 u^2 + \cdots$.*

(i) *Suppose that $F$ and $G$ are complex-valued measurable functions on $(0, \infty)$ such that $\int_0^\infty ds\, s^{\delta-1} e^{-\lambda_0 s} |F(s)| < \infty$ and $\int_0^\infty dz\, z^{-(\beta+1)} e^{-\lambda_0 z} |G(z)| < \infty$ for some $\lambda_0 > 0$. If $F$ and $G$ are connected with each other in such a way that*

$$\lambda^\delta \int_0^\infty ds\, s^{\delta-1} e^{-\lambda s} F(s) = R\left(\lambda^{-\beta} \int_0^\infty dz\, \frac{e^{-\lambda z}}{z^{\beta+1}} G(z)\right) \tag{3.1}$$

*for sufficiently large $\lambda > 0$, then $F(s)$ coincides with*

$$\sum_{n=1}^\infty \frac{r_n}{\Gamma(\delta+\beta n)} \int_{\Delta_n} \prod_{i=1}^n \frac{G(sv_i)}{v_i^{\beta+1}} \left(1 - \sum_{j=1}^n v_j\right)^{\delta+\beta n-1} dv_1 \cdots dv_n \qquad (\beta+\delta > 0) \tag{3.2}$$

*or*

$$r_1 G(s) + \sum_{n=2}^\infty r_n \int_{\Delta_{n-1}} \prod_{i=1}^{n-1} \frac{G(sv_i)}{v_i} \cdot \frac{G(s(1 - \sum_{j=1}^{n-1} v_j))}{1 - \sum_{j=1}^{n-1} v_j} dv_1 \cdots dv_{n-1}$$
$$(\beta = 0 = \delta) \tag{3.3}$$

*for a.e. $s > 0$, where the series and integrals converge absolutely for a.e. $s > 0$.*

(ii) *Suppose that $G$ is a complex-valued measurable function on $(0, \infty)$ such that $\int_0^\infty dz\, z^{-(\beta+1)} e^{-\lambda_0 z} |G(z)| < \infty$ for some $\lambda_0 > 0$. Then the series (3.2) or (3.3) with $|G(\cdot)|$ and $|r_n|$ in place of $G(\cdot)$ and $r_n$, respectively, converges for a.e. $s > 0$ and $F(s)$ defined by the expression (3.2) or (3.3) satisfies $\int_0^\infty ds\, s^{\delta-1} e^{-\lambda_0 s} |F(s)| < \infty$. Moreover, the relation (3.1) holds for sufficiently large $\lambda > 0$.*

**Proof.** We consider only the case where $\beta + \delta > 0$ since the other case can be handled quite similarly. To prove assertion (i), we first show that

$$G_0(s) := \sum_{n=1}^\infty \frac{|r_n|}{\Gamma(\delta+\beta n)} \int_{\Delta_n} \prod_{i=1}^n \frac{|G(sv_i)|}{v_i^{\beta+1}} \left(1 - \sum_{j=1}^n v_j\right)^{\delta+\beta n-1} dv_1 \cdots dv_n$$

is finite for a.e. $s > 0$. Setting, for each $s > 0$ and $n = 1, 2, \ldots,$

$$\Delta_n(s) = \{(z_1, \ldots, z_n) : z_1 \geq 0, \ldots, z_n \geq 0, z_1 + \cdots + z_n \leq s\}, \tag{3.4}$$

we observe, by Fubini's theorem, that

$$\int_0^\infty \mathrm{d}s \, s^{\delta-1} \mathrm{e}^{-\lambda s} \int_{\Delta_n} \prod_{i=1}^n \frac{|G(sv_i)|}{v_i^\beta} \left(1 - \sum_{j=1}^n v_j\right)^{\delta+\beta n-1} \frac{\mathrm{d}v_1 \cdots \mathrm{d}v_n}{v_1 \cdots v_n}$$

$$= \int_0^\infty \mathrm{d}s \, \mathrm{e}^{-\lambda s} \int_{\Delta_n(s)} \prod_{i=1}^n \frac{|G(z_i)|}{z_i^\beta} \left(s - \sum_{j=1}^n z_j\right)^{\delta+\beta n-1} \frac{\mathrm{d}z_1 \cdots \mathrm{d}z_n}{z_1 \cdots z_n}$$

$$= \int_{(0,\infty)^n} \mathrm{d}z_1 \cdots \mathrm{d}z_n \prod_{i=1}^n \frac{|G(z_i)|}{z_i^{\beta+1}} \int_{\sum_{j=1}^n z_j}^\infty \mathrm{d}s \, \mathrm{e}^{-\lambda s} \left(s - \sum_{j=1}^n z_j\right)^{\delta+\beta n-1}$$

$$= \frac{\Gamma(\delta+\beta n)}{\lambda^\delta} \left(\lambda^{-\beta} \int_0^\infty \mathrm{d}z \, \frac{\mathrm{e}^{-\lambda z}}{z^{\beta+1}} |G(z)|\right)^n.$$

Therefore, term-by-term integration yields

$$\lambda^\delta \int_0^\infty \mathrm{d}s \, s^{\delta-1} \mathrm{e}^{-\lambda s} G_0(s) = \sum_{n=1}^\infty |r_n| \left(\lambda^{-\beta} \int_0^\infty \mathrm{d}z \, \frac{\mathrm{e}^{-\lambda z}}{z^{\beta+1}} |G(z)|\right)^n.$$

Since the above series converges for sufficiently large $\lambda$, $G_0(s) < \infty$ for a.e. $s > 0$ and thus the desired absolute convergence follows.

Calculations needed for the proof of (i) have almost been done. Indeed, denoting by $G_1(s)$ the sum (3.2), one can see by obvious modification of the above calculations that (3.1) with $F$ replaced by $G_1$ holds for sufficiently large $\lambda > 0$. We thus conclude, by virtue of the uniqueness of Laplace transforms, that $F(s) = G_1(s)$ for a.e. $s > 0$. This proves assertion (i). The proof of assertion (ii) is essentially contained in the above and is therefore omitted. $\qquad\square$

At least formally, (3.3) is obtained as the 'degenerate limit' of (3.2), that is, by first setting $\beta = 0$ in (3.2) and then taking the limit as $\delta \downarrow 0$. Although our inversion formula given in Lemma 3.1(i) requires the Laplace transform to have a certain special form, one advantage is that it is described in the 'real world', that is, we do not need any complex integrals. Considered as prototypes of the inversion formula are integral representations of the Dickman function and the Buchstab function (cf. [1], page 22). Both functions appeared in asymptotic number theory [4, 12] and are related to PD(0, 1). On the other hand, assertion (ii) will be used below to compute the Laplace transform of a probability generating function of the PD($\alpha, \theta$) process discussed in [45].

In such applications, we will employ

$$R_{\alpha,\theta}(u) := \begin{cases} \Gamma(\theta)(e^{\theta u} - 1), & \alpha = 0, \, \theta > 0, \\ \Gamma(\theta + 1)\{(1 - C_\alpha u)^{-\theta/\alpha} - 1\}\theta^{-1}, & 0 < \alpha < 1, \, \theta \neq 0, \\ -\alpha^{-1}\log(1 - C_\alpha u), & 0 < \alpha < 1, \, \theta = 0, \end{cases} \tag{3.5}$$

where $C_\alpha = \alpha/\Gamma(1 - \alpha)$. The precise meanings of the power and logarithm in (3.5) as complex functions are as follows. Given $p \in \mathbf{R}$ and $t \in \mathbf{C} \setminus (-\infty, 0]$, define $t^p = \exp(p \log t)$, where $\log t = \log|t| + \arg t$, with $\arg t$ being chosen in $(-\pi, \pi)$. As mentioned in Exercise 1.2.7 of [44], the function $R_{\alpha,\theta}$ itself appears in connection with generalized Stirling numbers. Note that three expressions in the right-hand side of (3.5) vary continuously in $(\alpha, \theta)$, even in the limit as $\theta \to 0$ or $\alpha \downarrow 0$. This fact is consistent with the continuity of the two-parameter Poisson–Dirichlet family [52], [50]. More importantly, $R_{\alpha,\theta}$ admits an expansion of the form

$$R_{\alpha,\theta}(u) = \sum_{n=1}^\infty \frac{\Gamma(\theta + \alpha n)c_{n,\alpha,\theta}}{n!} u^n \tag{3.6}$$

as long as $|u| < 1/C_\alpha (= \infty$ for $\alpha = 0$ by definition). Combining Corollary 2.2 and Lemma 3.1, we prove the following formula, a refinement of Corollaries 49 and 50 in the aforementioned paper [45] by Pitman and Yor. We will use the convention that $\inf \varnothing = \infty$ so that, for a subset $X$ of $\mathbf{R}$, $\inf X < \infty$ means $X \neq \varnothing$.

**Theorem 3.2.** *Suppose that* $g : (0, \infty) \to \mathbf{C}$ *is a measurable function such that* $\lambda_\alpha(g) := \inf\{\lambda > 0 : \int_0^\infty dz\, z^{-(\alpha+1)} e^{-\lambda z}|g(z) - 1| < \infty\} < \infty$. *Put*

$$\lambda_\alpha^*(g) = \inf\left\{\lambda' > \lambda_\alpha(g) : \frac{C_\alpha}{\lambda^\alpha} \int_0^\infty dz\, \frac{e^{-\lambda z}}{z^{\alpha+1}}(g(z) - 1) \notin [1, \infty) \text{ for all } \lambda \geq \lambda'\right\}$$

*so that, in particular,* $\lambda_0^*(g) = \lambda_0(g)$. *Then for a.e.* $s > 0$,

$$E_{\alpha,\theta}\left[\prod_{i=1}^\infty (1 + |g(sV_i) - 1|)\right] < \infty \tag{3.7}$$

*and for* $\lambda > \lambda_\alpha^*(g)$,

$$\frac{\lambda^\theta}{\Gamma(\theta + 1)} \int_0^\infty ds\, s^{\theta-1} e^{-\lambda s}\left(E_{\alpha,\theta}\left[\prod_{i=1}^\infty g(sV_i)\right] - 1\right)$$

$$= \begin{cases} \dfrac{1}{\theta}\exp\left(\theta\int_0^\infty dz\,\dfrac{e^{-\lambda z}}{z}(g(z) - 1)\right) - \dfrac{1}{\theta}, & \alpha = 0, \, \theta > 0, \\[3mm] \dfrac{1}{\theta}\left(1 - \dfrac{C_\alpha}{\lambda^\alpha}\int_0^\infty dz\,\dfrac{e^{-\lambda z}}{z^{\alpha+1}}(g(z) - 1)\right)^{-\theta/\alpha} - \dfrac{1}{\theta}, & 0 < \alpha < 1, \, \theta \neq 0, \\[3mm] -\dfrac{1}{\alpha}\log\left(1 - \dfrac{C_\alpha}{\lambda^\alpha}\int_0^\infty dz\,\dfrac{e^{-\lambda z}}{z^{\alpha+1}}(g(z) - 1)\right), & 0 < \alpha < 1, \, \theta = 0. \end{cases} \tag{3.8}$$

**Proof.** Consider

$$F_0(s) := E_{\alpha,\theta}\left[\prod_{i=1}^{\infty}(1 + |g(sV_i) - 1|)\right] - 1,$$

which admits a series expansion due to (2.11). Noting (3.6), we apply the first half of Lemma 3.1(ii) with $\beta = \alpha$, $\delta = \theta$ and $R = R_{\alpha,\theta}$ to show that $F_0(s) < \infty$ for a.e. $s > 0$. Also, the last half of Lemma 3.1(ii) can be applied to the series expression of $F(s) := E_{\alpha,\theta}[\prod_{i=1}^{\infty} g(sV_i)] - 1$ due to (2.12) and we obtain

$$\lambda^{\theta}\int_0^{\infty}\mathrm{d}s\, s^{\theta-1}\mathrm{e}^{-\lambda s}F(s) = R_{\alpha,\theta}\left(\lambda^{-\alpha}\int_0^{\infty}\mathrm{d}z\,\frac{\mathrm{e}^{-\lambda z}}{z^{\alpha+1}}(g(z)-1)\right)$$

for sufficiently large $\lambda$. This extends to all $\lambda > \lambda^*_{\alpha}(g)$ by analytic continuation, showing (3.8). The proof of Theorem 3.2 is thus complete. □

This result will be exploited later in a variety of ways by taking

$$g(s) = \begin{cases} \mathbf{1}_{(0,1)}(s), & \text{in Section 4,} \\ \exp(-s^p), & \text{in Section 6.2,} \\ \psi_{\nu}(\pm s), & \text{in Section 7,} \end{cases}$$

where $p > 0$ and $\psi_{\nu}$ is the characteristic function of a distribution $\nu$. We conclude this section with the observation that the probability generating function of the $\mathrm{PD}(\alpha,\theta)$ process can also be characterized by an integral equation.

**Proposition 3.3.** *Suppose that* $g:(0,\infty) \to \mathbf{C}$ *is as in Theorem 3.2 and set* $\phi(s) = g(s) - 1$. *Then* $\Pi_g(s) := E_{\alpha,\theta}[\prod_{i=1}^{\infty} g(sV_i)]$ *solves the following integral equation for a.e.* $s > 0$:

(i) *for* $\alpha = 0$ *and* $\theta > 0$,

$$s^{\theta}\Pi_g(s) - \theta\int_0^s g(s-t)t^{\theta-1}\Pi_g(t)\,\mathrm{d}t = 0; \tag{3.9}$$

(ii) *for* $0 < \alpha < 1$ *and* $\theta \neq 0$,

$$\int_0^s \mathrm{d}t\,(s-t)^{-\alpha}t^{\theta}\Pi_g(t) - \int_0^s \mathrm{d}t\,(s-t)^{-\alpha}t^{\theta}g(s-t)$$
$$- \alpha\int_0^s \mathrm{d}u\int_0^u \mathrm{d}t\,(u-t)^{-(\alpha+1)}\phi(u-t)t^{\theta}\Pi_g(t)$$
$$- \theta\int_0^s \mathrm{d}u\int_0^u \mathrm{d}t\,(u-t)^{-\alpha}g(u-t)t^{\theta-1}(\Pi_g(t)-1) = 0; \tag{3.10}$$

(iii) *for $0 < \alpha < 1$ and $\theta = 0$,*

$$\int_0^s dt\,(s-t)^{-\alpha}\Pi_g(t) - \int_0^s t^{-\alpha}g(t)\,dt$$

$$- \alpha \int_0^s du \int_0^u dt\,(u-t)^{-(\alpha+1)}\phi(u-t)\Pi_g(t) = 0. \tag{3.11}$$

**Proof.** Given a measurable function $f$ on $(0,\infty)$, let $\widehat{f}$ denote the Laplace transform of $f$ and introduce the temporary notation $f_\beta(z) = f(z)/z^\beta$ for $\beta \in \mathbf{R}$. All equalities below involving $\lambda$ hold at least for sufficiently large $\lambda$.

(i) Set $f(s) = s^\theta \Pi_g(s)$ so that $f_1(s) = s^{\theta-1}\Pi_g(s)$. We start with a simplified version of (3.8) for $\alpha = 0$, that is, $\lambda^\theta \widehat{f}_1(\lambda) = \Gamma(\theta)\exp(\theta\widehat{\phi}_1(\lambda))$. Taking the logarithmic derivative of both sides with respect to $\lambda$, we get

$$\frac{\theta}{\lambda} - \frac{\widehat{f}(\lambda)}{\widehat{f}_1(\lambda)} = -\theta\widehat{\phi}(\lambda)\left(-\theta\widehat{g}(\lambda) + \frac{\theta}{\lambda}\right)$$

and hence $\widehat{f}(\lambda) = \theta\widehat{g}(\lambda)\widehat{f}_1(\lambda)$. This suffices to prove (3.9).

(ii) For notational simplicity, let $\Phi(s) = \Pi_g(s) - 1$ and $G(s) = \Pi_g(s)$. In the case where $0 < \alpha < 1$ and $\theta \neq 0$, the equality (3.8) reads

$$\frac{\theta\lambda^\theta}{\Gamma(1+\theta)}\widehat{\Phi}_{1-\theta}(\lambda) = \left(1 - \frac{C_\alpha}{\lambda^\alpha}\widehat{\phi}_{\alpha+1}(\lambda)\right)^{-\theta/\alpha} - 1.$$

Setting $\beta = -\theta(<\alpha < 1)$, we convert the above identity to

$$(\Gamma(1-\beta)\lambda^\beta - \beta\widehat{\Phi}_{\beta+1}(\lambda))^{1/\beta} = (\Gamma(1-\alpha)\lambda^\alpha - \alpha\widehat{\phi}_{\alpha+1}(\lambda))^{1/\alpha}.$$

By taking the logarithmic derivative,

$$\frac{\widehat{\Phi}_\beta(\lambda) + \Gamma(1-\beta)\lambda^{\beta-1}}{\Gamma(1-\beta)\lambda^\beta - \beta\widehat{\Phi}_{\beta+1}(\beta)} = \frac{\widehat{\phi}_\alpha(\lambda) + \Gamma(1-\alpha)\lambda^{\alpha-1}}{\Gamma(1-\alpha)\lambda^\alpha - \alpha\widehat{\phi}_{\alpha+1}(\lambda)}.$$

Here, considering the function $\mathbf{1}_\gamma(z) := 1/z^\gamma$ for every $\gamma < 1$, note that the numerator of the left-hand side equals $\widehat{\Phi}_\beta(\lambda) + \widehat{\mathbf{1}}_\beta(\lambda) = \widehat{G}_\beta(\lambda)$ and, similarly,

$$\widehat{\phi}_\alpha(\lambda) + \Gamma(1-\alpha)\lambda^{\alpha-1} = \widehat{g}_\alpha(\lambda). \tag{3.12}$$

Consequently, we have

$$\widehat{\mathbf{1}}_\alpha(\lambda)\widehat{G}_\beta(\lambda) - \widehat{\mathbf{1}}_\beta(\lambda)\widehat{g}_\alpha(\lambda) - \frac{\alpha}{\lambda}\widehat{\phi}_{\alpha+1}(\lambda)\widehat{G}_\beta(\lambda) + \frac{\beta}{\lambda}\widehat{g}_\alpha(\lambda)\widehat{\Phi}_{\beta+1}(\lambda) = 0. \tag{3.13}$$

This proves (3.10).

(iii) Differentiating (3.8) with $\theta = 0$ and again using (3.12), we get

$$\widehat{\mathbf{1}}_{\alpha}(\lambda)\widehat{\Pi}_g(\lambda) - \frac{1}{\lambda}\widehat{g}_{\alpha}(\lambda) - \frac{\alpha}{\lambda}\widehat{\phi}_{\alpha+1}(\lambda)\widehat{\Pi}_g(\lambda) = 0, \tag{3.14}$$

which shows (3.11). The proof is thus completed. □

In fact, the generality of the assumptions in Proposition 3.3 makes the resulting equations (3.10) and (3.11) rather complicated. Under an additional hypothesis which ensures that both $g$ and $\Pi_g$ are of bounded variation on each finite interval, the following equations are derived instead.

For $0 < \alpha < 1$ and $\theta \neq 0$,

$$\int_0^s (s-t)^{-\alpha}t^{\theta}\,\mathrm{d}\Pi_g(t) - \int_0^s (s-t)^{\theta}t^{-\alpha}\,\mathrm{d}g(t)$$
$$- \int_0^s \mathrm{d}t\,\{(s-t)^{-\alpha}\theta t^{\theta-1} + \alpha(s-t)^{-(\alpha+1)}t^{\theta}\}\phi(s-t)(\Pi_g(t) - 1) = 0. \tag{3.15}$$

For $0 < \alpha < 1$ and $\theta = 0$,

$$\int_0^s (s-t)^{-\alpha}\,\mathrm{d}\Pi_g(t) - s^{-\alpha}\phi(s) - \alpha\int_0^s \mathrm{d}t\,(s-t)^{-(\alpha+1)}\phi(s-t)\Pi_g(t) = 0. \tag{3.16}$$

The proof changes only after multiplying (3.13) and (3.14) by $\lambda$ and then uses

$$\widehat{\mathbf{1}}_{\alpha}(\lambda)\lambda\widehat{G}_{\beta}(\lambda) - \widehat{\mathbf{1}}_{\beta}(\lambda)\lambda\widehat{g}_{\alpha}(\lambda) = \widehat{\mathbf{1}}_{\alpha}(\lambda)\int_0^{\infty} \mathrm{e}^{-\lambda s}\,\mathrm{d}\Phi_{\beta}(s) - \widehat{\mathbf{1}}_{\beta}(\lambda)\int_0^{\infty} \mathrm{e}^{-\lambda s}\,\mathrm{d}\phi_{\alpha}(s)$$

and $\lambda\widehat{\Pi}_g(\lambda) = 1 + \int_0^{\infty} \mathrm{e}^{-\lambda s}\,\mathrm{d}\Pi_g(s)$, respectively. The details are omitted.

# 4. The two-parameter generalization of the Dickman function

One of fundamental 'observables' in a point process on $\mathbf{R}$ is the 'position of the last particle' (cf. [27]), if any. As for the $\mathrm{PD}(\alpha, \theta)$ process, this is nothing but $V_1$, the *first* component of a $\mathrm{PD}(\alpha, \theta)$-distributed random element $(V_i)$. For special values of $\alpha$ and $\theta$, the law of $V_1$ was found in various contexts much earlier than Kingman's discovery of the Poisson–Dirichlet limit. For example, the Dickman function [12], usually denoted $\rho(\cdot)$, is identified with $\rho(s) = P_{0,1}(sV_1 < 1)$; see, for example, Section III 5.3 of [48] for related discussions in asymptotic number theory and Section 1.1 of [1]. Here and in what follows, $P_{\alpha,\theta}$ is a probability distribution under which $(V_i)$ is $\mathrm{PD}(\alpha, \theta)$-distributed. Also, a distribution function found independently in [8] and [35] can be identified with $P_{\alpha,0}(V_1^{-1} \leq s)$. See also [32] ((77), page 14). It is natural to introduce the two-parameter version of the Dickman function by

$$\rho_{\alpha,\theta}(s) = P_{\alpha,\theta}(sV_1 < 1).$$

Clearly, $\rho_{\alpha,\theta}(s) = 1$ for all $s \leq 1$. The one-parameter family $\{\rho_{0,\theta} : \theta > 0\}$ has been studied in connection with both population genetics [18, 54] and the asymptotic theory of the symmetric group [22]. For details, we refer the reader to the identity (4.13) in [3] (resp., Lemma 4.7 in [1]), where $g_\theta(s)$ (resp., $p_\theta(s)$) is identical with $s^{\theta-1}\rho_{0,\theta}(s)$ up to some multiplicative constant. It also appears in a natural extension [21] of Dickman's result in number theory. Another interesting context in which the one-parameter family arises is the identification of limit distributions associated with random minimal directed spanning trees [40].

The aim of this section is to describe consequences for $\rho_{\alpha,\theta}$'s which follow from the results in Sections 2 and 3. First, a choice of $\phi$ in (2.12) to give an expression for $\rho_{\alpha,\theta}$ is $\phi(v) = -\mathbf{1}_{[1,\infty)}(sv)$, with $s > 0$ being given. With this choice, (2.12) now reads

$$\rho_{\alpha,\theta}(s) = \sum_{n=0}^{\infty} \frac{(-1)^n c_{n,\alpha,\theta}}{n!} I_{n,\alpha,\theta}(s), \tag{4.1}$$

where $I_{0,\alpha,\theta}(\cdot) \equiv 1$ and where, for $n = 1, 2, \ldots,$

$$I_{n,\alpha,\theta}(s) = \int_{\Delta_n} \prod_{i=1}^{n} \frac{\mathbf{1}_{[1,\infty)}(sv_i)}{v_i^{\alpha+1}} \left(1 - \sum_{j=1}^{n} v_j\right)^{\theta+\alpha n - 1} \mathrm{d}v_1 \cdots \mathrm{d}v_n. \tag{4.2}$$

The reader is cautioned that our notation $I_{n,\alpha,\theta}$ is in conflict with that of [45]. We observe that $I_{n,\alpha,\theta}(s) = 0$ whenever $n > s$. So, the right-hand side of (4.1) is in fact a finite sum taken over $n$ with $0 \leq n \leq s$. In the case $(\alpha, \theta) = (0, 1)$, we recover a well-known formula for the Dickman function (cf. (1.35) in [1]), while the above identity with $\theta = 0$ extends the formula for $P_{\alpha,0}(V_1^{-1} \leq s) =: H_\alpha(s)$, obtained in [35] ((3.7), page 730), to all values of $s \geq 1$. As a two-parameter example, we find in [34] (Theorem 3.3.1, page 164) this type of expression for $\rho_{1/2,1/2}$, which determines the limit distribution of the maximal size of trees associated with a random mapping.

It is proved in Propositions 19 and 20 of [45] that for all $s > 0$,

$$\rho_{\alpha,\theta}(s) = \frac{\Gamma(\theta+1)}{\Gamma(\theta+\alpha)\Gamma(1-\alpha)} \int_0^{\min\{s^{-1},1\}} \mathrm{d}v \, \frac{(1-v)^{\theta+\alpha-1}}{v^{\alpha+1}} \rho_{\alpha,\theta+\alpha}\left(\frac{1-v}{v}\right) \tag{4.3}$$

and hence the probability density $P_{\alpha,\theta}(V_1 \in \mathrm{d}v)/\mathrm{d}v$ is given by

$$f_{\alpha,\theta}(v) := \frac{\Gamma(\theta+1)}{\Gamma(\theta+\alpha)\Gamma(1-\alpha)} \cdot \frac{(1-v)^{\theta+\alpha-1}}{v^{\alpha+1}} \rho_{\alpha,\theta+\alpha}\left(\frac{1-v}{v}\right) \mathbf{1}_{(0,1)}(v). \tag{4.4}$$

For later convenience, we remark that underlying (4.3) are 'termwise equalities' shown in the next lemma. For every positive integer $n$ and $s > 0$, let

$$\nabla_n(s) = \{(v_1, \ldots, v_n) : v_1 \geq \cdots \geq v_n \geq 0, v_1 + \cdots + v_n \leq s\}.$$

In particular, $\nabla_n(1)$ is simply denoted by $\nabla_n$.

**Lemma 4.1.** *For all $n = 1, 2, \ldots$ and $s > 0$, it holds that*

$$I_{n,\alpha,\theta}(s) = n \int_{\min\{s^{-1}, 1\}}^{1} \mathrm{d}v \, \frac{(1-v)^{\theta+\alpha-1}}{v^{\alpha+1}} I_{n-1,\alpha,\theta+\alpha}\left(\frac{1-v}{v}\right). \tag{4.5}$$

**Proof.** Obviously, we may assume that $s \geq 1$. Also, (4.5) with $n = 1$ is clear from the definition (4.2). For $n \geq 2$, by symmetry of the integrand in (4.2), $I_{n,\alpha,\theta}(s)$ equals

$$n! \int_0^1 \frac{\mathrm{d}v_n}{v_n^{\alpha+1}} \mathbf{1}_{\{sv_n \geq 1\}} \int_{\nabla_{n-1}(1-v_n)} \frac{\mathrm{d}v_1 \cdots \mathrm{d}v_{n-1}}{v_1^{\alpha+1} \cdots v_{n-1}^{\alpha+1}} \mathbf{1}_{\{v_{n-1} \geq v_n\}} \left(1 - \sum_{j=1}^{n} v_j\right)^{\theta+\alpha n-1}. \tag{4.6}$$

The change of variables $v_i = (1 - v_n)u_i$ converts the inner integral to

$$(1-v_n)^{\theta+\alpha-1} \int_{\nabla_{n-1}} \frac{\mathrm{d}u_1 \cdots \mathrm{d}u_{n-1}}{u_1^{\alpha+1} \cdots u_{n-1}^{\alpha+1}} \mathbf{1}_{\{(1-v_n)u_{n-1} \geq v_n\}} \left(1 - \sum_{j=1}^{n-1} u_j\right)^{\theta+\alpha n-1}$$

$$= \frac{1}{(n-1)!}(1-v_n)^{\theta+\alpha-1} I_{n-1,\alpha,\theta+\alpha}\left(\frac{1-v_n}{v_n}\right).$$

Substituting this into (4.6), we get (4.5). $\qquad\square$

It should be noted that (4.3) does not give a closed equation for $\rho_{\alpha,\theta}$ unless $\alpha = 0$. To derive such equations for the general case, we shall apply some results from the previous section. Setting $g = \mathbf{1}_{(0,1)}$ in (3.8) and observing that $\lambda_\alpha^*(g) = 0$, we see from Theorem 3.2 that for all $\lambda > 0$,

$$\lambda^\theta \int_0^\infty \mathrm{d}s \, s^{\theta-1} \mathrm{e}^{-\lambda s}(\rho_{\alpha,\theta}(s) - 1) = R_{\alpha,\theta}\left(-\lambda^{-\alpha} \int_1^\infty \mathrm{d}z \, \frac{\mathrm{e}^{-\lambda z}}{z^{\alpha+1}}\right). \tag{4.7}$$

For $\theta > 0$, this relation is rewritten in a slightly simpler form:

$$\frac{\lambda^\theta}{\Gamma(\theta)} \int_0^\infty \mathrm{d}s \, s^{\theta-1} \mathrm{e}^{-\lambda s} \rho_{\alpha,\theta}(s) = \begin{cases} \exp\left(-\theta \int_1^\infty \mathrm{d}z \, \frac{\mathrm{e}^{-\lambda z}}{z}\right), & \alpha = 0, \\ \left(1 + \frac{C_\alpha}{\lambda^\alpha} \int_1^\infty \mathrm{d}z \, \frac{\mathrm{e}^{-\lambda z}}{z^{\alpha+1}}\right)^{-\theta/\alpha}, & 0 < \alpha < 1. \end{cases} \tag{4.8}$$

Integral equations satisfied by $\rho_{\alpha,\theta}$ will be derived as the ones equivalent to (4.7), in other words, as consequences of Proposition 3.3, its variant (3.15) and (3.16). Preliminary observations concerning this are that $\Pi_g(s) = \rho_{\alpha,\theta}(s)$ for $g = \mathbf{1}_{(0,1)}$ and that (3.9), (3.15) and (3.16) hold for all $s \geq 1$ because of continuity. As will be discussed, the equation below is known in the cases $\alpha = 0$ and $\theta = 0$. For this reason, the proof of the next theorem will be concerned with the remaining case only.

**Theorem 4.2.** $\rho_{\alpha,\theta}$ *solves the equation*

$$(s-1)^\theta \rho_{\alpha,\theta}(s-1) + \int_{s-1}^s (s-t)^{-\alpha} t^\theta \, \mathrm{d}\rho_{\alpha,\theta}(t) = 0 \qquad (s > 1). \tag{4.9}$$

**Proof.** The case where $0 < \alpha < 1$ and $\theta \neq 0$. (3.15) now reads, for all $s > 1$, as

$$\int_0^s (s-t)^{-\alpha} t^\theta \, \mathrm{d}\rho_{\alpha,\theta}(t) + (s-1)^\theta$$

$$+ \int_0^{s-1} \mathrm{d}t \, \{(s-t)^{-\alpha} \theta t^{\theta-1} + \alpha(s-t)^{-(\alpha+1)} t^\theta\} (\rho_{\alpha,\theta}(t) - 1) = 0$$

because $-\mathrm{d}g(t)$ on $(0, \infty)$ equals $\delta_1(\mathrm{d}t)$. Integration by parts yields (4.9). $\qquad\square$

*Remarks.* (i) (4.9) with $\alpha = 0$ or its variant can be found in [19, 21, 22, 54]; see also (4.25) in [1].

(ii) One can easily see that (4.9) with $\theta = 0$ is derived from a functional equation for $H_\alpha(t) = 1 - \rho_{\alpha,0}(t)$ obtained in [35] ((3.5), page 730).

In the next section, we calculate not only the marginal distributions $P_{\alpha,\theta}(V_m \in \cdot)$ for $m = 2, 3, \ldots$, but also multidimensional distributions of PD$(\alpha,\theta)$, by developing point process calculus based on Theorem 2.1.

# 5. Distributional results for PD$(\alpha, \theta)$

In this section, we apply the theory of point processes to deduce from (2.3) some distributional information on PD$(\alpha,\theta)$. An essential idea underlying the subsequent argument is the principle of inclusion-exclusion, which was already being used in [19]. By the same reasoning as in the proof of Theorem 2 there, namely, by a version of Fréchet's formula (see, for example, Section IV.5 in [15]), the following relationship between one-dimensional distributions and correlation functions holds with great generality.

**Lemma 5.1.** *Let* $-\infty \leq a < b \leq \infty$ *and let* $\xi$ *be a simple point process on* $(a, b)$ *with correlation measures* $\mu_1, \mu_2, \ldots$ *. Suppose that* $\xi((a, b)) = \infty$ *a.s. and that* $\xi([c, b)) < \infty$ *a.s. for each* $c \in (a, b)$. *Given a positive integer* $m$, *let* $Z_m$ *be the* $m$th *largest point in* $\xi$. *Then, for every* $z \in (a, b)$,

$$P(Z_m \geq z) = \frac{1}{(m-1)!} \sum_{k=0}^\infty \frac{(-1)^k}{(m+k)k!} \int_{[z,b)^{m+k}} \mu_{m+k}(\mathrm{d}y_1 \cdots \mathrm{d}y_{m+k}). \tag{5.1}$$

*If, in addition,* $\xi$ *has correlation functions* $q_1, q_2, \ldots$, *then*

$$\frac{P(Z_m \in \mathrm{d}z)}{\mathrm{d}z} = \frac{\mathbf{1}_{(a,b)}(z)}{(m-1)!} \sum_{k=0}^\infty \frac{(-1)^k}{k!} \int_{[z,b)^{m+k-1}} \mathrm{d}y_1 \cdots \mathrm{d}y_{m+k-1} \, q_{m+k}(z, y_1, \ldots, y_{m+k-1}).$$

**Example 5.1.** Consider a Poisson point process $Z_1 > Z_2 > \cdots$ on $(a, b)$ with mean measure $\Lambda(\mathrm{d}z) = h(z)\,\mathrm{d}z$ such that $\Lambda((a,b)) = \infty$ and $\Lambda([c,b)) < \infty$ for all $c \in (a, b)$. Lemma 5.1 then gives

$$\frac{P(Z_m \in \mathrm{d}z)}{\mathrm{d}z} = \frac{h(z)}{(m-1)!}\left(\int_z^b h(y)\,\mathrm{d}y\right)^{m-1}\exp\left(-\int_z^b h(y)\,\mathrm{d}y\right)\mathbf{1}_{(a,b)}(z) \tag{5.2}$$

because the $n$th correlation function is $h(y_1)\cdots h(y_n)$; see, for example, Example 2.5 in [27].

Setting $\rho_{m,\alpha,\theta}(s) = P_{\alpha,\theta}(sV_m < 1)$ $(m = 1, 2, \ldots)$, we obtain the following proposition containing a two-parameter generalization of the aforementioned result in [19].

**Proposition 5.2.** *Let $m$ be a positive integer. Then, for all $s > 0$,*

$$\rho_{m,\alpha,\theta}(s) = 1 - \frac{1}{(m-1)!}\sum_{0 \le k \le s-m}\frac{(-1)^k}{(m+k)k!}c_{m+k,\alpha,\theta}I_{m+k,\alpha,\theta}(s) \tag{5.3}$$

*and for all $\lambda > 0$,*

$$\lambda^\theta\int_0^\infty \mathrm{d}s\, s^{\theta-1}\mathrm{e}^{-\lambda s}(\rho_{m,\alpha,\theta}(s) - 1) = R_{m,\alpha,\theta}\left(-\lambda^{-\alpha}\int_1^\infty \mathrm{d}z\,\frac{\mathrm{e}^{-\lambda z}}{z^{\alpha+1}}\right), \tag{5.4}$$

*where*

$$R_{m,\alpha,\theta}(u) := \frac{(-1)^{m-1}c_{m-1,\alpha,\theta}}{(m-1)!}\int_0^u x^{m-1}R_{\alpha,\theta+\alpha(m-1)}{}'(x)\,\mathrm{d}x$$

$$= \begin{cases} \dfrac{\Gamma(\theta+1)(-\theta)^{m-1}}{(m-1)!}\displaystyle\int_0^u x^{m-1}\mathrm{e}^{\theta x}\,\mathrm{d}x, & \alpha = 0, \\[2ex] \dfrac{\Gamma(\theta+1)\Gamma(\theta/\alpha+m)(-\alpha)^{m-1}}{(m-1)!\Gamma(\theta/\alpha+1)\Gamma(1-\alpha)^m}\displaystyle\int_0^u \frac{x^{m-1}}{(1-C_\alpha x)^{\theta/\alpha+m}}\,\mathrm{d}x, & 0 < \alpha < 1. \end{cases} \tag{5.5}$$

*Also, for $v \in (0, 1)$,*

$$\frac{P_{\alpha,\theta}(V_m \in \mathrm{d}v)}{\mathrm{d}v} = \frac{\Gamma(\theta+1)}{\Gamma(\theta+\alpha)\Gamma(1-\alpha)}v^{-(\alpha+1)}(1-v)^{\theta+\alpha-1}$$

$$\times \frac{1}{(m-1)!}\sum_{0 \le k \le 1/v-m}\frac{(-1)^k}{k!}c_{m+k-1,\alpha,\theta+\alpha}I_{m+k-1,\alpha,\theta+\alpha}\left(\frac{1-v}{v}\right). \tag{5.6}$$

**Remarks.** (i) Since $\rho_{1,\alpha,\theta} = \rho_{\alpha,\theta}$, it is worth noting that the right-hand sides of (5.3), (5.5) and (5.6) with $m = 1$ equal the right-hand sides of (3.5), (4.1) and (4.4), respectively.

(ii) In [34], Theorem 3 on page 47 gives the expression (5.3) with $(\alpha, \theta) = (0, 1)$ to the limit distribution of the $m$th maximal cycle length in a random permutation.

(iii) At the end of this section, another expression of the density (5.6) will be given in terms of the generalized Dickman function.

**Proof of Proposition 5.2.** First, (5.1) and Theorem 2.1 together imply that

$$P_{\alpha,\theta}(sV_m \geq 1) = \frac{1}{(m-1)!} \sum_{k=0}^{\infty} \frac{(-1)^k}{(m+k)k!} c_{m+k,\alpha,\theta} I_{m+k,\alpha,\theta}(s).$$

By removing the terms which actually vanish, (5.3) follows. Next, showing (5.4) reduces to verifying that for $u$ sufficiently close to 0,

$$R_{m,\alpha,\theta}(u) = \frac{(-1)^{m-1}}{(m-1)!} \sum_{k=0}^{\infty} \frac{\Gamma(\theta + \alpha(m+k))}{(m+k)k!} c_{m+k,\alpha,\theta} u^{m+k}. \tag{5.7}$$

For, if (5.7) is true, then (5.3) and Lemma 3.1(ii) with $R = R_{m,\alpha,\theta}$ and $G = -\mathbf{1}_{[1,\infty)}$ show that (5.4) is valid, at least for $\lambda$ sufficiently large, and the expression (5.5) allows us to extend (5.4) to all $\lambda > 0$ by analytic continuation. Turning to (5.7), one can easily verify it by using (term-by-term differentiation of) (3.6) or, alternatively, by substituting the expansion of $e^{\theta x}$ or $(1 - C_\alpha x)^{-(\theta/\alpha + m)}$ into (5.5). This proves (5.4).

Lastly, since (5.3) is at hand, the proof of (5.6) is similar to that of (4.4) based on (4.5). The details are left to the reader. □

To obtain multidimensional results for $PD(\alpha,\theta)$, one needs more developed arguments still based essentially on inclusion-exclusion. In the theory of point processes, one such calculus is formulated as a connection between correlation measures and Jonassy measures; see Section 5 of [7]. Its significance is that the local probabilistic structure of points in the process is revealed in terms of correlation measures.

**Lemma 5.3.** *Let* $-\infty \leq a < b \leq \infty$ *and let* $\xi$ *be a simple point process on* $(a,b)$ *with correlation functions* $q_1, q_2, \ldots$. *Suppose that* $\xi((a,b)) = \infty$ *a.s. and that for each* $c \in (a,b)$, *the nth factorial moments* $M_n(c)$ *of* $\xi([c,b))$ *satisfy*

$$\sum_{n=1}^{\infty} \frac{M_n(c)}{n!} (1+\varepsilon)^n < \infty \qquad \text{for some } \varepsilon = \varepsilon(c) > 0. \tag{5.8}$$

*Let* $Z_1 > Z_2 > \cdots$ *be the decreasing sequence of points in* $\xi$. *Then, for each* $m = 1, 2, \ldots$, *the joint probability density* $f_m$ *of* $(Z_1, \ldots, Z_m)$ *is given by*

$$f_m(\mathbf{z}) = \sum_{k=0}^{\infty} \frac{(-1)^k}{k!} \int_{[z_m,b)^k} dy_1 \cdots dy_k \, q_{m+k}(\mathbf{z}, y_1, \ldots, y_k), \tag{5.9}$$

*where* $\mathbf{z} = (z_1, \ldots, z_m)$ *is such that* $b > z_1 > \cdots > z_m > a$.

**Proof.** Let $b =: z_0 > z_1 > \cdots > z_m > a$ be given arbitrarily. Consider the localized process $\xi^{(m)}(\mathrm{d}x) := \xi(\mathrm{d}x \cap [z_m, b))$, the $n$th correlation function of which is given by $q_n^{(m)}(\mathbf{x}_n) := q_n(\mathbf{x}_n)\mathbf{1}_{[z_m, b)}(x_1) \cdots \mathbf{1}_{[z_m, b)}(x_n)$, where $\mathbf{x}_n = (x_1, \ldots, x_n)$. By applying Theorem 5.4.II in [7], then, we have

$$P(Z_1 > z_1 \geq Z_2 > z_2 \geq \cdots > z_{m-2} \geq Z_{m-1} > z_{m-1}, z_m > Z_m)$$

$$= P(\xi^{(m)}([z_m, z_0]) = m - 1, \xi^{(m)}((z_1, z_0]) = \cdots = \xi^{(m)}((z_{m-1}, z_{m-2}]) = 1)$$

$$= \sum_{k=0}^{\infty} \frac{(-1)^k}{k!} \int_{(a,b)^{m-1+k}} \prod_{i=1}^{m-1} \mathbf{1}_{(z_i, z_{i-1})}(x_i) q_{m-1+k}^{(m)}(\mathbf{x}_{m-1+k}) \, \mathrm{d}x_1 \cdots \mathrm{d}x_{m-1+k}$$

$$= \int_{z_{m-1}}^{z_{m-2}} \mathrm{d}x_{m-1} \cdots \int_{z_2}^{z_1} \mathrm{d}x_2 \int_{z_1}^{b} \mathrm{d}x_1$$

$$\times \sum_{k=0}^{\infty} \frac{(-1)^k}{k!} \int_{[z_m, b)^k} \mathrm{d}x_m \, \mathrm{d}x_{m+1} \cdots \mathrm{d}x_{m-1+k} \, q_{m-1+k}(\mathbf{x}_{m-1+k}) \tag{5.10}$$

$$= \int_{z_{m-1}}^{z_{m-2}} \mathrm{d}x_{m-1} \cdots \int_{z_2}^{z_1} \mathrm{d}x_2 \int_{z_1}^{b} \mathrm{d}x_1$$

$$\times \sum_{k=0}^{\infty} (-1)^k \int_{z_m}^{b} \mathrm{d}x_m \int_{x_m}^{b} \mathrm{d}x_{m+1} \cdots \int_{x_{m-2+k}}^{b} \mathrm{d}x_{m-1+k} \, q_{m-1+k}(\mathbf{x}_{m-1+k}), \tag{5.11}$$

where the last equality is due to symmetry in $x_m, \ldots, x_{m-1+k}$ and the terms for $k = 0$ in (5.10) and (5.11) are understood as $q_{m-1}(\mathbf{x}_{m-1})$ with the convention that $q_0 \equiv 1$. This shows (5.9) by the symmetry of the integrand in (5.11) in $k-1$ variables $x_{m+1}, \ldots, x_{m-1+k}$. $\qquad\square$

*Example 5.2.* If $\xi$ is a Poisson point process with mean measure density $h$ as in Example 5.1, the joint density of the first $m$ largest points in $\xi$ is given by

$$\sum_{k=0}^{\infty} \frac{(-1)^k}{k!} \int_{[z_m, b)^k} \mathrm{d}y_1 \cdots \mathrm{d}y_k \prod_{i=1}^{m} h(z_i) \prod_{j=1}^{k} h(y_j) = \prod_{i=1}^{m} h(z_i) \exp\left(-\int_{z_m}^{b} h(z) \, \mathrm{d}z\right), \tag{5.12}$$

provided that $b > z_1 > \cdots > z_m > a$.

By applying (5.9), we deduce the next result for joint densities, which is essentially contained in [45], although an explicit formula is not given. (Indeed, the authors obtained in their Corollary 41 a corresponding result for the variables $y_1, y_2, \ldots$ in (2.6) instead of $v_1, v_2, \ldots$. Alternatively, the formula below can be retrieved from their Proposition 47 after some additional calculations. See also Lemma 3.1 in [16].) This simultaneously generalizes the formula due to Watterson [54] for $\alpha = 0$ and the one-dimensional result (4.4). Also, our expression (5.13) below will be quite useful in Section 6.1.

**Theorem 5.4.** *Let $(V_i)$ be governed by* $\mathrm{PD}(\alpha, \theta)$*. Then, for each $m = 1, 2, \ldots$, the joint probability density of $(V_1, \ldots, V_m)$ at $\mathbf{v} = (v_1, \ldots, v_m) \in \nabla_m$ is*

$$f_{m,\alpha,\theta}(\mathbf{v}) := c_{m,\alpha,\theta} \prod_{i=1}^{m} v_i^{-(\alpha+1)} \left(1 - \sum_{j=1}^{m} v_j\right)^{\theta + \alpha m - 1} \rho_{\alpha, \theta + \alpha m}\left(\frac{1 - \sum_{j=1}^{m} v_j}{v_m}\right). \quad (5.13)$$

**Proof.** First note that the $\mathrm{PD}(\alpha, \theta)$ process $\xi = \sum \delta_{V_i}$ is a point process on $(0, 1)$ which satisfies all of the assumptions in Lemma 5.3 since, for each $0 < c < 1$, $\xi([c, 1)) \leq 1/c$ a.s. by $V_i \leq 1/i$ $(i = 1, 2, \ldots)$. According to (5.9), the density to be computed is

$$\sum_{k=0}^{\infty} \frac{(-1)^k}{k!} \int_{[v_m, 1)^k} \mathrm{d}u_1 \cdots \mathrm{d}u_k \, q_{m+k, \alpha, \theta}(\mathbf{v}, u_1, \ldots, u_k).$$

Here, by (2.3) and (2.2),

$$\int_{[v_m, 1)^k} \mathrm{d}u_1 \cdots \mathrm{d}u_k \, q_{m+k, \alpha, \theta}(\mathbf{v}, u_1, \ldots, u_k)$$

$$= c_{m, \alpha, \theta} \prod_{i=1}^{m} v_i^{-(\alpha+1)} c_{k, \alpha, \theta + \alpha m} \quad (5.14)$$

$$\times \int_{\Delta_k (1 - \sum_{i=1}^{m} v_i)} \prod_{j=1}^{k} \frac{\mathbf{1}_{\{u_j \geq v_m\}}}{u_j^{\alpha+1}} \left(1 - \sum_{i=1}^{m} v_i - \sum_{j=1}^{k} u_j\right)^{\theta + \alpha(m+k) - 1} \mathrm{d}u_1 \cdots \mathrm{d}u_k$$

with notation (3.4). Because the last integral in (5.14) is equal to

$$\left(1 - \sum_{i=1}^{m} v_i\right)^{\theta + \alpha m - 1} I_{k, \alpha, \theta + \alpha m}\left(\frac{1 - \sum_{i=1}^{m} v_i}{v_m}\right), \quad (5.15)$$

the desired expression (5.13) is derived from (5.14), (5.15) and (4.1). $\qquad \square$

Since $f_{m, \alpha, \theta}$ is too complicated to compute quantities concerning the joint distribution directly from this, it is worth providing a moment formula. For the one-parameter case where $\alpha = 0$ and $\theta > 0$, Griffiths [18] showed that $E_{0, \theta}[V_1^{a_1} \cdots V_m^{a_m}]$ equals

$$\frac{\theta^m \Gamma(\theta)}{\Gamma(\theta + a_1 + \cdots + a_m)} \int_{z_1 > \cdots > z_m > 0} \mathrm{d}\mathbf{z} \prod_{i=1}^{m} (z_i^{a_i - 1} \mathrm{e}^{-z_i}) \exp\left(-\theta \int_{z_m}^{\infty} \mathrm{d}z \, \frac{\mathrm{e}^{-z}}{z}\right), \quad (5.16)$$

provided that $\theta + a_1 + \cdots + a_m > 0$, where $\mathrm{d}\mathbf{z} = \mathrm{d}z_1 \cdots \mathrm{d}z_m$. The proof reduces to calculus of a gamma process by virtue of the well-known independence property of it; see [50] and [51] for extensive discussions and related topics. In this regard, we remark only that the integrand in (5.16) contains the density function, computed from (5.12), of the first $m$ largest points in a gamma process with parameter $\theta$, that is, a Poisson point process on $(0, \infty)$ with mean measure $\theta \, dz/(ze^z)$. Although the case where $0 < \alpha < 1$ and $\theta > -\alpha$ can be handled by means of representing a $\mathrm{PD}(\alpha, \theta)$-distributed random element by an $\alpha$-stable subordinator ([45, 50] combined with (5.12)), we prefer to exploit the previous results in order to make the proof self-contained.

**Corollary 5.5.** *Let* $0 < \alpha < 1$ *and* $\theta > -\alpha$. *If* $a_1 + \cdots + a_m > -\theta$, *then*

$$
\begin{aligned}
&E_{\alpha, \theta}[V_1{}^{a_1} \cdots V_m{}^{a_m}] \\
&= \frac{\Gamma(\theta + 1)\Gamma(\theta/\alpha + m)\alpha^{m-1}}{\Gamma(\theta + a_1 + \cdots + a_m)\Gamma(\theta/\alpha + 1)\Gamma(1 - \alpha)^m} \\
&\quad \times \int_{z_1 > \cdots > z_m > 0} d\mathbf{z} \prod_{i=1}^m (z_i^{a_i - (\alpha+1)} e^{-z_i}) \left(1 + C_\alpha \int_{z_m}^\infty dz \, \frac{e^{-z}}{z^{\alpha+1}}\right)^{-(\theta/\alpha+m)}.
\end{aligned}
\tag{5.17}
$$

**Proof.** By (5.13),

$$
\begin{aligned}
&E_{\alpha, \theta}[V_1{}^{a_1} \cdots V_m{}^{a_m}]/c_{m,\alpha,\theta} \\
&= \int_{\nabla_m} dv_1 \cdots dv_m \prod_{i=1}^m v_i^{a_i - (\alpha+1)} \left(1 - \sum_{j=1}^m v_j\right)^{\theta + \alpha m - 1} \rho_{\alpha, \theta + \alpha m}\left(\frac{1 - \sum_{j=1}^m v_j}{v_m}\right).
\end{aligned}
$$

Multiplying both sides by $\Gamma(\theta + A) = \int_0^\infty y^{\theta + A - 1} e^{-y} \, dy$ with $A = a_1 + \cdots + a_m$, we introduce an additional integration with respect to a new variable $y$ on the right-hand side. For the resulting $(m+1)$-dimensional integral, perform the change of variables

$$
s := v_m^{-1}\left(1 - \sum_{j=1}^m v_j\right), \qquad z_i := v_i y \qquad (i = 1, \ldots, m),
$$

or, equivalently, $y = (z_1 + \cdots + z_m) + z_m s, v_i = z_i/y \ (i = 1, \ldots, m)$, to get

$$
\begin{aligned}
&E_{\alpha, \theta}[V_1{}^{a_1} \cdots V_m{}^{a_m}]\Gamma(\theta + A)/c_{m,\alpha,\theta} \\
&= \int_{z_1 > \cdots > z_m > 0} d\mathbf{z} \prod_{i=1}^m z_i^{a_i - (\alpha+1)} z_m^{\theta + \alpha m - 1} \int_0^\infty ds \, y^m e^{-y} |J| s^{\theta + \alpha m - 1} \rho_{\alpha, \theta + \alpha m}(s)
\end{aligned}
$$

with Jacobian $J = \partial(v_1, \ldots, v_m, y)/\partial(z_1, \ldots, z_m, s) = y^{-m} z_m$. Thanks to (4.8), with $\theta + \alpha m > 0$ in place of $\theta$, we arrive at (5.17). $\qquad \square$

**Remarks.** (i) It can be seen that the formula for $E_{\alpha, \theta}[V_m^p] (p > -\theta)$ in Proposition 17 of [45] is recovered by putting $a_1 = \cdots = a_{m-1} = 0$ and $a_m = p$ in (5.17).

(ii) Let $0 < v < 1$ and $m \in \{2, 3, \ldots\}$. By integrating $f_{m,\alpha,\theta}(v_1, \ldots, v_{m-1}, v)$ with respect to $dv_1 \cdots dv_{m-1}$ over $\{1 > v_1 > \cdots > v_{m-1} > v\}$, we can deduce a representation (an alternative to (5.6)) of the density $P_{\alpha,\theta}(V_m \in \mathrm{d}v)/\mathrm{d}v$ of the form

$$\frac{c_{m,\alpha,\theta}(1-v)^{\theta+\alpha-1}}{(m-1)!v^{\alpha+1}}\int_{\Delta_{m-1}}\mathrm{d}u_1\cdots\mathrm{d}u_{m-1}\prod_{i=1}^{m-1}\frac{\mathbf{1}_{\{u_i\geq v/(1-v)\}}}{u_i^{\alpha+1}}u_m^{\theta+\alpha m-1}\rho_{\alpha,\theta+\alpha m}\left(\frac{1-v}{v}u_m\right),$$

where $u_m = 1 - (u_1 + \cdots + u_{m-1})$. This formula, a natural extension of (4.4), is verified by calculations similar to (5.14) and (5.15). The details are left to the reader.

# 6. Asymptotics of $\mathrm{PD}(\alpha,\theta)$ for large $\theta$

In this section, we study certain asymptotic behaviors of $\mathrm{PD}(\alpha,\theta)$ as $\theta \to \infty$, generalizing results of Griffiths [18] and of Joyce, Krone and Kurtz [28], who all worked on $\mathrm{PD}(0,\theta)$'s with motivation coming from the study of population genetics. Although there is a context [24] in which such an extension could be applicable, this section is mainly intended to demonstrate that the results we have thus far obtained provide efficient methods for the study of two-parameter Poisson–Dirichlet distributions. It turns out that the presence of the parameter $\alpha$ does not affect the validity of the assertions analogous to those for $\mathrm{PD}(0,\theta)$'s, except for some minor changes of sub-leading terms or multiplicative constants in rescaling.

## 6.1. Convergence to a Gumbel point process

Given $0 \leq \alpha < 1$ and $\theta > 1$, put

$$\beta_{\alpha,\theta} = \log\theta - (\alpha+1)\log\log\theta - \log\Gamma(1-\alpha). \tag{6.1}$$

Under the assumption that $(V_i^{(\theta)})_{i=1}^\infty$ is distributed according to $\mathrm{PD}(0,\theta)$, it was shown in [18] that as $\theta \to \infty$, $(\theta V_i^{(\theta)} - \beta_{0,\theta})_{i=1}^\infty$ converges in law to a Poisson point process with mean measure density $\exp(-z)$, $-\infty < z < \infty$, which may be called a *Gumbel point process* since the largest point in the process obeys the Gumbel distribution $\exp(-\exp(-z))$, as seen from (5.2); see also (6.3) below. A two-parameter generalization of this result is the following.

**Theorem 6.1.** *Fix an $\alpha \in [0,1)$ arbitrarily. For each $\theta > 1$, let $(V_i^{(\alpha,\theta)})_{i=1}^\infty$ have $\mathrm{PD}(\alpha,\theta)$ distribution and define*

$$Z_i^{(\alpha,\theta)} = \theta V_i^{(\alpha,\theta)} - \beta_{\alpha,\theta} \qquad (i = 1, 2, \ldots). \tag{6.2}$$

*Then, as $\theta \to \infty$, $(Z_1^{(\alpha,\theta)}, Z_2^{(\alpha,\theta)}, \ldots)$ converges in joint distribution to the decreasing order statistics $(Z_1^*, Z_2^*, \ldots)$ of a Poisson point process on the whole real axis with mean measure $\exp(-z)\,\mathrm{d}z$.*

The proof of this theorem reduces to showing pointwise convergence of joint densities, thanks to Scheffé's theorem (cf. [2]). It follows from (5.12) that

$$P(Z_1^* \in \mathrm{d}z_1, \ldots, Z_m^* \in \mathrm{d}z_m) = \exp\left(-\sum_{i=1}^m z_i - \mathrm{e}^{-z_m}\right) \mathrm{d}z_1 \cdots \mathrm{d}z_m, \qquad (6.3)$$

while the density of $(Z_1^{(\alpha,\theta)}, \ldots, Z_m^{(\alpha,\theta)})$ derived from (5.13) involves the generalized Dickman function as a key factor. It is therefore essential to study the asymptotic behavior of $\rho_{\alpha,\theta}$ with a properly rescaled variable for large $\theta$. This is the content of the following lemma, in which is derived the Gumbel distribution, the law of $Z_1^*$.

**Lemma 6.2.** *Fix $0 \le \alpha < 1$. Let $\theta_0 > 1$ and $b \colon [\theta_0, \infty) \to \mathbf{R}$ be bounded. If $V_1^{(\alpha,\theta)}$ and $Z_1^{(\alpha,\theta)}$ are as in Theorem 6.1, then for each $x \in \mathbf{R}$,*

$$P(Z_1^{(\alpha,\theta)} < x - (x + \beta_{\alpha,\theta})b(\theta)V_1^{(\alpha,\theta)}) = \rho_{\alpha,\theta}\left(\frac{\theta}{x + \beta_{\alpha,\theta}} + b(\theta)\right) \to \mathrm{e}^{-\mathrm{e}^{-x}} \qquad (6.4)$$

*as $\theta \to \infty$. In particular, $Z_1^{(\alpha,\theta)} \to Z_1^*$ weakly and $\rho_{\alpha,\theta}(s) \to 1$ for all $s \in \mathbf{R}$.*

**Proof.** First, the equality in (6.4) is clear from (6.2). Rather than using (4.1), which is not very informative for the estimation of the value of $\rho_{\alpha,\theta}$, we employ (4.8), that is,

$$\lambda^\theta \int_0^\infty \mathrm{d}s\, s^{\theta-1} \mathrm{e}^{-\lambda s} \rho_{\alpha,\theta}(s) = \begin{cases} \Gamma(\theta) \exp\left(-\theta \int_\lambda^\infty \mathrm{d}z\, \dfrac{\mathrm{e}^{-z}}{z}\right), & \alpha = 0, \\ \Gamma(\theta)\left(1 + C_\alpha \int_\lambda^\infty \mathrm{d}z\, \dfrac{\mathrm{e}^{-z}}{z^{\alpha+1}}\right)^{-\theta/\alpha}, & 0 < \alpha < 1, \end{cases} \qquad (6.5)$$

where $\lambda > 0$ is arbitrary. We will give only a proof for the case where $0 < \alpha < 1$ because the case $\alpha = 0$ can be handled in much the same way. By Fubini's theorem, the left-hand side of (6.5) is written as

$$\lambda^\theta E_{\alpha,\theta}\left[\int_0^{V_1^{-1}} \mathrm{d}s\, s^{\theta-1} \mathrm{e}^{-\lambda s}\right] = E_{\alpha,\theta}\left[\int_0^{\lambda V_1^{-1}} \mathrm{d}t\, t^{\theta-1} \mathrm{e}^{-t}\right] = E_{\alpha,\theta}[f_{\lambda,\theta}(V_1^{-1})],$$

where $f_{\lambda,\theta}(s) = \int_0^{\lambda s} \mathrm{d}t\, t^{\theta-1} \mathrm{e}^{-t}$. For all $\lambda > 0$ and $s > 0$, we get

$$\rho_{\alpha,\theta}(s) \le \frac{E_{\alpha,\theta}[f_{\lambda,\theta}(V_1^{-1})]}{f_{\lambda,\theta}(s)} = \frac{\Gamma(\theta)}{f_{\lambda,\theta}(s)}\left(1 + C_\alpha \int_\lambda^\infty \mathrm{d}z\, \frac{\mathrm{e}^{-z}}{z^{\alpha+1}}\right)^{-\theta/\alpha} \qquad (6.6)$$

as a Chebyshev-type bound. Set $\varepsilon(\theta) = \theta^{-1/4}$, so that

$$\varepsilon(\theta)\log\theta \to 0, \qquad \varepsilon(\theta)\sqrt{\theta} \to \infty \qquad (\theta \to \infty). \qquad (6.7)$$

Given arbitrary $x \in \mathbf{R}$, choose

$$s = s(x, \theta) := \frac{\theta}{x + \beta_{\alpha, \theta}} + b(\theta) \quad \text{and} \quad \lambda = \lambda(x, \theta) := (1 + \varepsilon(\theta))(x + \beta_{\alpha, \theta}),$$

both of which are positive for sufficiently large $\theta$. Then $\widetilde{\varepsilon}(x, \theta)$, defined implicitly by $\lambda(x, \theta)s(x, \theta) = (1 + \widetilde{\varepsilon}(x, \theta))\theta$, has the same properties as (6.7). The latter property, combined with the standard argument in the proof of Stirling's formula (that is, the Laplace method [9]), shows that

$$\frac{\Gamma(\theta)}{f_{\lambda(x,\theta),\theta}(s(x,\theta))} = \frac{\int_0^\infty \mathrm{d}t\, t^{\theta-1}\mathrm{e}^{-t}}{\int_0^{(1+\widetilde{\varepsilon}(x,\theta))\theta} \mathrm{d}t\, t^{\theta-1}\mathrm{e}^{-t}} = \frac{\int_0^\infty \mathrm{d}u\, u^{\theta-1}\mathrm{e}^{-\theta u}}{\int_0^{1+\widetilde{\varepsilon}(x,\theta)} \mathrm{d}u\, u^{\theta-1}\mathrm{e}^{-\theta u}} \to 1$$

as $\theta \to \infty$. Consequently, by (6.6), with the above choice of $s$ and $\lambda$, we have

$$\limsup_{\theta \to \infty} \rho_{\alpha,\theta}\left(\frac{\theta}{x + \beta_{\alpha,\theta}} + b(\theta)\right) \le \limsup_{\theta \to \infty}\left(1 + C_\alpha \int_{\lambda(x,\theta)}^\infty \mathrm{d}z\, \frac{\mathrm{e}^{-z}}{z^{\alpha+1}}\right)^{-\theta/\alpha}. \tag{6.8}$$

For two functions $a_1(\theta)$ and $a_2(\theta)$, which may have parameters $\alpha$, $x$, etcetera, the notation $a_1(\theta) \sim a_2(\theta)$ will mean that $a_1(\theta)/a_2(\theta) \to 1$ as $\theta \to \infty$. By (6.7) and (6.1),

$$\int_{\lambda(x,\theta)}^\infty \mathrm{d}z\, \frac{\mathrm{e}^{-z}}{z^{\alpha+1}} \sim \frac{\mathrm{e}^{-\lambda(x,\theta)}}{\lambda(x,\theta)^{\alpha+1}} \sim \frac{\mathrm{e}^{-(x+\beta_{\alpha,\theta})}}{(x+\beta_{\alpha,\theta})^{\alpha+1}} \sim \frac{\Gamma(1-\alpha)}{\theta}\mathrm{e}^{-x}.$$

Combining this with (6.8) yields

$$\limsup_{\theta \to \infty} \rho_{\alpha,\theta}\left(\frac{\theta}{x + \beta_{\alpha,\theta}} + b(\theta)\right) \le \exp(-\mathrm{e}^{-x}). \tag{6.9}$$

The converse estimate can be shown in an analogous way. Indeed, considering this time $1 - \rho_{\alpha,\theta}(s) = P_{\alpha,\theta}(V_1 \ge s^{-1})$, we again have, by Fubini's theorem,

$$\lambda^\theta \int_0^\infty \mathrm{d}s\, s^{\theta-1}\mathrm{e}^{-\lambda s}(1 - \rho_{\alpha,\theta}(s)) = E_{\alpha,\theta}[\varphi_{\lambda,\theta}(V_1)], \tag{6.10}$$

where $\varphi_{\lambda,\theta}(s) = \int_{\lambda/s}^\infty \mathrm{d}t\, t^{\theta-1}\mathrm{e}^{-t}$. Therefore, a Chebyshev-type bound we get is

$$1 - \rho_{\alpha,\theta}(s) \le \frac{E_{\alpha,\theta}[\varphi_{\lambda,\theta}(V_1)]}{\varphi_{\lambda,\theta}(s^{-1})} = \frac{\Gamma(\theta)}{\int_{\lambda s}^\infty \mathrm{d}t\, t^{\theta-1}\mathrm{e}^{-t}}\left\{1 - \left(1 + C_\alpha \int_\lambda^\infty \mathrm{d}z\, \frac{\mathrm{e}^{-z}}{z^{\alpha+1}}\right)^{-\theta/\alpha}\right\},$$

where the last equality follows from (6.10) and (6.5) together. By setting $s = s(x, \theta)$ and $\lambda = (1 - \varepsilon(\theta))(x + \beta_{\alpha,\theta})$, one can easily modify the previous argument to obtain

$$\limsup_{\theta \to \infty}\left\{1 - \rho_{\alpha,\theta}\left(\frac{\theta}{x + \beta_{\alpha,\theta}} + b(\theta)\right)\right\} \le 1 - \exp(-\mathrm{e}^{-x}).$$

This, together with (6.9), proves the convergence in (6.4). The other assertions can be seen simply by taking $b(\cdot) \equiv 0$ and noting that $\rho_{\alpha,\theta}(\cdot)$ is non-increasing. $\qquad \square$

**Proof of Theorem 6.1.** Fix a positive integer $m$ and $z_1 > \cdots > z_m$ arbitrarily. By (5.13), the joint probability density of $(Z_1^{(\alpha,\theta)}, \ldots, Z_m^{(\alpha,\theta)})$ at $(z_1, \ldots, z_m)$ is

$$
\theta^{-m} f_{m,\alpha,\theta}\left( \frac{z_1 + \beta_{\alpha,\theta}}{\theta}, \ldots, \frac{z_m + \beta_{\alpha,\theta}}{\theta} \right)
$$

$$
= c_{m,\alpha,\theta} \prod_{i=1}^{m} (z_i + \beta_{\alpha,\theta})^{-(\alpha+1)} \left( 1 - \sum_{j=1}^{m} \frac{z_j + \beta_{\alpha,\theta}}{\theta} \right)^{\theta + \alpha m - 1} \tag{6.11}
$$

$$
\times \theta^{\alpha m} \rho_{\alpha,\theta+\alpha m}\left( \frac{\theta - \sum_{j=1}^{m}(z_j + \beta_{\alpha,\theta})}{z_m + \beta_{\alpha,\theta}} \right) \mathbf{1}_{\nabla_m}\left( \frac{z_1 + \beta_{\alpha,\theta}}{\theta}, \ldots, \frac{z_m + \beta_{\alpha,\theta}}{\theta} \right).
$$

It is easy to see that the factors in (6.11) behave as

$$
c_{m,\alpha,\theta} \prod_{i=1}^{m} (z_i + \beta_{\alpha,\theta})^{-(\alpha+1)} \sim \theta^{(1-\alpha)m} \Gamma(1-\alpha)^{-m} (\log \theta)^{-(1+\alpha)m}
$$

and

$$
\left( 1 - \sum_{j=1}^{m} \frac{z_j + \beta_{\alpha,\theta}}{\theta} \right)^{\theta + \alpha m - 1} \sim \exp\left( -\sum_{j=1}^{m} z_j \right) \theta^{-m} (\log \theta)^{(1+\alpha)m} \Gamma(1-\alpha)^m,
$$

respectively. In view of (6.3), it remains to show that the factor involving $\rho_{\alpha,\theta+\alpha m}(\cdots)$ converges to $\exp(-\exp(-z_m))$. But this follows from (6.4) with $b(\cdot)$ defined appropriately. The proof of Theorem 6.1 is thus complete. $\qquad \square$

The weak convergence of $Z_m^{(\alpha,\theta)}$ to $Z_m^*$ implies the following.

**Corollary 6.3.** *For each $m = 1, 2, \ldots$ and $x \in \mathbf{R}$, as $\theta \to \infty$,*

$$
\rho_{m,\alpha,\theta}\left( \frac{\theta}{x + \beta_{\alpha,\theta}} \right) \to \exp(-e^{-x}) \sum_{k=0}^{m-1} \frac{e^{-kx}}{k!}. \tag{6.12}
$$

**Proof.** As a consequence of Theorem 6.1, the left-hand side of (6.12), that is, $P(Z_m^{(\alpha,\theta)} \leq x)$, tends to $P(Z_m^* \leq x)$, which is calculated by (5.2):

$$
P(Z_m^* \leq x) = \int_{-\infty}^{x} \frac{\mathrm{d}y}{(m-1)!} \exp(-my - e^{-y}) = \frac{1}{(m-1)!} \int_{e^{-x}}^{\infty} t^{m-1} e^{-t} \, \mathrm{d}t.
$$

This changes into the right-hand side of (6.12) after repeated integration by parts. $\qquad \square$

Theorem 6.1 immediately implies that $(V_1^{(\alpha,\theta)}, V_2^{(\alpha,\theta)}, \ldots) \to (0, 0, \ldots)$ as $\theta \to \infty$. Feng [16] has recently obtained associated large deviation estimates with rate function $(v_i) \mapsto -\log(1 - \sum v_i)$. The proof is based on logarithmic asymptotics of the finite-dimensional densities. We note that, in view of (5.13), one can grasp the validity of such a result because the factor involving $\rho_{\alpha, \theta + \alpha m}$ stays away from 0 by (6.4).

## 6.2. Asymptotic normality of the population moments

As the final topic on $\mathrm{PD}(\alpha, \theta)$ for large $\theta$, we discuss asymptotic normality of the population moment, which is defined for each $p > 0$ by the sum

$$H_p^{(\alpha,\theta)} = (V_1^{(\alpha,\theta)})^p + (V_2^{(\alpha,\theta)})^p + \cdots,$$

where $(V_1^{(\alpha,\theta)}, V_2^{(\alpha,\theta)}, \ldots)$ continues to be $\mathrm{PD}(\alpha, \theta)$-distributed. We can verify that $H_p^{(\alpha,\theta)} < \infty$, a.s. for all $p > \alpha$. Indeed, setting $g_p(s) = \exp(-s^p)$, one observes that $\lambda_\alpha(g_p)$ in Theorem 3.2 is finite for $p > \alpha$ and the claim follows from (3.7) with $g = g_p$. Also, (2.12) with $\phi(v) = g_p(sv) - 1$ makes it possible to derive the Laplace transform of the law of $H_p^{(\alpha,\theta)}$ because $\prod_{i=1}^{\infty} g_p(sV_i^{(\alpha,\theta)}) = \exp(-s^p H_p^{(\alpha,\theta)})$. As far as the $\mathrm{PD}(0, \theta)$ case is concerned, the statistics $H_2^{(0,\theta)}$, called the *population homozygosity*, is known to play a special role in some population genetics contexts (see [19, 20, 29]), where $V_i^{(0,\theta)}$ are regarded as the ranked frequencies of alleles in a population. More generally, $H_p^{(0,\theta)}$ is referred to as the *pth population moment* in [28], where the asymptotic normality of $H_p^{(0,\theta)}$ for large $\theta$ has been established. On the other hand, covariances of $H_p^{(\alpha,0)}$'s are calculated in [10] ((12), page 188), motivated by certain physical problems. In this subsection, we will provide an extension of these results to the two-parameter setting.

In the case $\alpha = 0$, it was shown in [28] that for each $p = 2, 3, \ldots$, as $\theta \to \infty$, the limit distribution of $W_p^{(\theta)} := \sqrt{\theta}(\theta^{p-1} H_p^{(0,\theta)} / \Gamma(p) - 1)$ is the normal distribution with mean 0 and variance $\Gamma(2p)/\Gamma(p)^2 - p^2$. In fact, in [28], much stronger results were obtained, showing, for example, asymptotic normality of random vectors $(W_m^{(\theta)})_{m=2}^{\infty}$. To seek an appropriate rescaling for the two-parameter case, we give some auxiliary results by means of the first and second correlation functions, namely, $q_{1,\alpha,\theta}$ and $q_{2,\alpha,\theta}$.

**Lemma 6.4.** *Let $p > \alpha$ and $p' > \alpha$. Then*

$$E[H_p^{(\alpha,\theta)}] = \frac{\Gamma(\theta+1)\Gamma(p-\alpha)}{\Gamma(\theta+p)\Gamma(1-\alpha)} =: h_p^{(\alpha,\theta)} \tag{6.13}$$

*and the covariance* $\mathrm{Cov}(H_p^{(\alpha,\theta)}, H_{p'}^{(\alpha,\theta)})$ *of $H_p^{(\alpha,\theta)}$ and $H_{p'}^{(\alpha,\theta)}$ is given by*

$$h_{p+p'}^{(\alpha,\theta)} + \frac{\Gamma(p-\alpha)\Gamma(p'-\alpha)}{\Gamma(1-\alpha)^2}\left(\frac{\Gamma(\theta+1)(\theta+\alpha)}{\Gamma(\theta+p+p')} - \frac{\Gamma(\theta+1)^2}{\Gamma(\theta+p)\Gamma(\theta+p')}\right). \tag{6.14}$$

*Moreover, as $\theta \to \infty$,*

$$\theta^{p-1} E[H_p^{(\alpha,\theta)}] \to \frac{\Gamma(p-\alpha)}{\Gamma(1-\alpha)} \tag{6.15}$$

*and*

$$\theta^{p+p'-1} \operatorname{Cov}(H_p^{(\alpha,\theta)}, H_{p'}^{(\alpha,\theta)}) \to \frac{\Gamma(p+p'-\alpha)}{\Gamma(1-\alpha)} + \frac{\Gamma(p-\alpha)\Gamma(p'-\alpha)}{\Gamma(1-\alpha)^2}(\alpha - pp'). \tag{6.16}$$

**Proof.** (6.13) and (6.14) are immediate by calculating $\int_0^1 v^p q_{1,\alpha,\theta}(v) \, dv$ and

$$E[H_p^{(\alpha,\theta)} H_{p'}^{(\alpha,\theta)}] = \int_0^1 v^{p+p'} q_{1,\alpha,\theta}(v) \, dv + \int_{\Delta_2} v_1^p v_2^{p'} q_{2,\alpha,\theta}(v_1, v_2) \, dv_1 \, dv_2,$$

*respectively. By virtue of the well-known formula $\Gamma(\theta+1)/\Gamma(\theta+p) \sim \theta^{1-p}$ implied by Stirling's formula, (6.15) follows from (6.13). Showing (6.16), however, requires a more accurate version of Stirling's formula: with the standard o-notation, $\Gamma(\theta+1) = \theta^\theta \sqrt{2\pi\theta} \exp(-\theta + C\theta^{-1} + o(\theta^{-1}))$ as $\theta \to \infty$ for some universal constant $C$; see, for example, [9]. This makes it possible to show that*

$$\frac{\Gamma(\theta+1)}{\Gamma(\theta+p)} = \theta^{1-p} \exp(-p(p-1)(2\theta)^{-1} + o(\theta^{-1})). \tag{6.17}$$

*The proof of (6.16) is a matter of straightforward calculation, combining (6.17) with (6.14). The details are omitted.* □

Lemma 6.4 suggests a rescaling of the form

$$W_p^{(\alpha,\theta)} := \sqrt{\theta} \left( \frac{\Gamma(1-\alpha)}{\Gamma(p-\alpha)} \theta^{p-1} H_p^{(\alpha,\theta)} - 1 \right). \tag{6.18}$$

We can now state a desired result for $0 < \alpha < 1$.

**Theorem 6.5.** *Fix an $\alpha \in (0,1)$ arbitrarily. Suppose that $p > \alpha$ and $p \neq 1$. For each $\theta > 0$, let $W_p^{(\alpha,\theta)}$ be defined as above. Then, as $\theta \to \infty$, $W_p^{(\alpha,\theta)}$ converges in law to a normal random variable with mean 0 and variance*

$$\sigma_{\alpha,p}^2 := \frac{\Gamma(1-\alpha)\Gamma(2p-\alpha)}{\Gamma(p-\alpha)^2} + \alpha - p^2 > 0. \tag{6.19}$$

**Proof.** *Fixing an arbitrary $x \in \mathbf{R}$, we will prove that as $\theta \to \infty$,*

$$\psi_p^{(\alpha,\theta)}(x) := E[\exp(\sqrt{-1}x W_p^{(\alpha,\theta)})] \to \exp\left(-\frac{\sigma_{\alpha,p}^2}{2} x^2\right). \tag{6.20}$$

The proof is based on (3.8), a consequence of Theorem 3.2 which can be written as

$$\frac{1}{\Gamma(\theta)} \int_0^\infty \mathrm{d}s\, s^{\theta-1} \mathrm{e}^{-s} E_{\alpha,\theta}\left[\prod_{i=1}^\infty g(sV_i)\right] = \left(1 + C_\alpha \int_0^\infty \mathrm{d}z\, \frac{\mathrm{e}^{-z}}{z^{\alpha+1}}(1 - g(z))\right)^{-\theta/\alpha}, \quad (6.21)$$

provided that $1 > \lambda_\alpha^*(g)$. Let $g$ be of the form $g(z) = \exp(c_1 z^p + d_1 z)$ and observe that

$$E_{\alpha,\theta}\left[\prod_{i=1}^\infty g(sV_i)\right] = E[\exp(c_1 s^p H_p^{(\alpha,\theta)} + d_1 s)]. \quad (6.22)$$

Note that (6.21) is valid as long as the integral on the right-hand side is small enough in modulus. With the choices

$$c_1 = \sqrt{-1}\,\frac{x}{\sqrt{\theta}} \cdot \frac{\Gamma(1-\alpha)}{\Gamma(p-\alpha)} \quad \text{and} \quad d_1 = -\sqrt{-1}\,\frac{x}{\sqrt{\theta}}, \quad (6.23)$$

this requirement for the validity of (6.21) is fulfilled for $\theta$ large enough since as $\theta \to \infty$,

$$\int_0^\infty \mathrm{d}z\, \frac{\mathrm{e}^{-z}}{z^{\alpha+1}}(1 - g(z)) = -\sqrt{-1}\,\frac{x}{\sqrt{\theta}} \int_0^\infty \mathrm{d}z\, \frac{\mathrm{e}^{-z}}{z^{\alpha+1}}\left(\frac{\Gamma(1-\alpha)}{\Gamma(p-\alpha)} z^p - z\right)$$

$$+ \frac{x^2}{2\theta} \int_0^\infty \mathrm{d}z\, \frac{\mathrm{e}^{-z}}{z^{\alpha+1}}\left(\frac{\Gamma(1-\alpha)}{\Gamma(p-\alpha)} z^p - z\right)^2 + \mathrm{o}(\theta^{-1})$$

$$= \frac{x^2}{2\theta}\Gamma(1-\alpha)(\sigma_{\alpha,p}^2 + (p-1)^2) + \mathrm{o}(\theta^{-1}),$$

where the last equality is achieved by direct calculations. This also implies that the right-hand side of (6.21) tends to $\exp(-\sigma_{\alpha,p}^2 x^2/2 - (p-1)^2 x^2/2)$. On the other hand, the left-hand side of (6.21) with (6.23) becomes, by (6.22), after change of variable $s = \theta t$,

$$\int_0^\infty \mathrm{d}t\, t^{\theta-1} \mathrm{e}^{-\theta t} E[\exp(c_1 \theta^p H_p^{(\alpha,\theta)} t^p + d_1 \theta t)](\Gamma(\theta)\theta^{-\theta})^{-1}$$

$$= \int_0^\infty \mathrm{d}t\, t^{-1} \mathrm{e}^{-\theta h(t)} \psi_p^{(\alpha,\theta)}(x t^p) \exp(\sqrt{-1}x\sqrt{\theta}(t^p - t))(\Gamma(\theta)\theta^{-\theta}\mathrm{e}^\theta)^{-1},$$

where $h(t) = t - \log t - 1$. Thus, by Stirling's formula, (6.20) is equivalent to

$$\frac{\int_0^\infty \mathrm{d}t\, t^{-1} \mathrm{e}^{-\theta h(t)} \psi_p^{(\alpha,\theta)}(x t^p) \exp(\sqrt{-1}x\sqrt{\theta}(t^p - t))}{\sqrt{2\pi/\theta}} - \psi_p^{(\alpha,\theta)}(x)\mathrm{e}^{-(p-1)^2 x^2/2} \to 0,$$

the proof of which will be divided into proofs of

$$\sqrt{\theta} \int_0^\infty \mathrm{d}t\, t^{-1} \mathrm{e}^{-\theta h(t)}(\psi_p^{(\alpha,\theta)}(x t^p) - \psi_p^{(\alpha,\theta)}(x)) \exp(\sqrt{-1}x\sqrt{\theta}(t^p - t)) \to 0 \quad (6.24)$$

and

$$\sqrt{\frac{\theta}{2\pi}} \int_0^\infty \mathrm{d}t \, t^{-1} \mathrm{e}^{-\theta h(t)} \exp(\sqrt{-1}x\sqrt{\theta}(t^p - t)) - \mathrm{e}^{-(p-1)^2 x^2/2} \to 0. \qquad (6.25)$$

Note that the integral in (6.24) can be bounded in modulus by

$$\int_0^\infty \mathrm{d}t \, t^{-1} \mathrm{e}^{-\theta h(t)} |\psi_p^{(\alpha,\theta)}(xt^p) - \psi_p^{(\alpha,\theta)}(x)| \le |x| E[|W_p^{(\alpha,\theta)}|] \int_0^\infty \mathrm{d}t \, t^{-1} \mathrm{e}^{-\theta h(t)} |t^p - 1|$$

and that $\sup_{\theta > 1} E[|W_p^{(\alpha,\theta)}|] < \infty$ since, by Lemma 6.4 and (6.17),

$$(E[|W_p^{(\alpha,\theta)}|])^2 \le \operatorname{Var}(W_p^{(\alpha,\theta)}) + (E[W_p^{(\alpha,\theta)}])^2 \to \sigma_{\alpha,p}^2.$$

Therefore, our task is reduced to establishing that $\sqrt{\theta} \int_0^\infty \mathrm{d}t \, t^{-1} \mathrm{e}^{-\theta h(t)} |t^p - 1| \to 0$ and (6.25). Both are purely analytic problems which can be solved without difficulty by the Laplace method (or by the change of variable $t = 1 + u/\sqrt{\theta}$).

It remains to show that $\sigma_{\alpha,p} > 0$. A stronger result which can be shown is

$$\sigma_{\alpha,p}{}^2 = (1-\alpha) \operatorname{Var}\left(\frac{\Gamma(1-\alpha)}{\Gamma(p-\alpha)} Y_\alpha^{p-1}\right) + \frac{\alpha}{1-\alpha}(p-1)^2, \qquad (6.26)$$

where $Y_\alpha$ has the density $z^{1-\alpha}\mathrm{e}^{-z}/\Gamma(2-\alpha), 0 < z < \infty$. The proof of (6.26) is fairly direct and therefore omitted. The proof of Theorem 6.5 is thus complete. $\qquad \square$

**Remark.** The above proof also works in the case $\alpha = 0$ if (6.21) is replaced by the corresponding equality for $\alpha = 0$. This gives a proof which is completely different from that in [28] and which is not based on the independence property underlying $\mathrm{PD}(0, \theta)$.

One may consider an extension of Theorem 1 in [28] in the two-parameter setting, that is, convergence of $\{W_p^{(\alpha,\theta)} : p > \alpha\}$ to a centered Gaussian system as $\theta \to \infty$. The convergence of every finite-dimensional law is shown by a suitable modification of the proof of Theorem 6.5. It means, for instance, that $g$ in (6.21) is taken to be

$$g(s) = \exp\left[\frac{\sqrt{-1}}{\sqrt{\theta}} \sum_{k=1}^n x_k\left(\frac{\Gamma(1-\alpha)}{\Gamma(p_k - \alpha)} s^{p_k} - s\right)\right],$$

where $x_k \in \mathbf{R}$, $p_1 > \cdots > p_n > \alpha$ and $n = 1, 2, \ldots$ are arbitrary. The limit process is described by its covariance $\Gamma(1-\alpha)\Gamma(p + p' - \alpha)/(\Gamma(p-\alpha)\Gamma(p'-\alpha)) + \alpha - pp' =: C(p, p')$ and the associated quadratic form $\sum C(p_k, p_l) x_k x_l$ has two different expressions:

$$\frac{1}{\Gamma(1-\alpha)} \int_0^\infty \mathrm{d}z \, \frac{\mathrm{e}^{-z}}{z^{\alpha+1}} \left\{\sum_{k=1}^n x_k\left(\frac{\Gamma(1-\alpha)}{\Gamma(p_k - \alpha)} z^{p_k} - z\right)\right\}^2 - \left\{\sum_{k=1}^n x_k(p_k - 1)\right\}^2$$

$$= (1-\alpha) \operatorname{Var}\left(\sum_{k=1}^n x_k \frac{\Gamma(1-\alpha)}{\Gamma(p_k - \alpha)} Y_\alpha^{p_k - 1}\right) + \frac{\alpha}{1-\alpha}\left\{\sum_{k=1}^n x_k(p_k - 1)\right\}^2,$$

where $Y_\alpha$ is as in (6.26). We do not produce the whole proof, a routine matter which is left to the reader.

# 7. The generalized Dirichlet process

This section discusses the two-parameter generalization of Dirichlet processes from the point of view of our previous results. There are a number of motivations to study the original Dirichlet process, for example, as a prior distribution in Bayesian nonparametric statistics [17] and as a stationary state of a certain diffusion process arising in population genetics. The well-known relationship between the Dirichlet process and the Poisson–Dirichlet distribution is described as follows. Let $\{X_i\}_{i=1}^\infty$ be a sequence of i.i.d. random variables and $(V_i)_{i=1}^\infty$ be governed by $\mathrm{PD}(0, \theta)$. Assuming independence of $\{X_i\}_{i=1}^\infty$ and $(V_i)_{i=1}^\infty$, we know that a random distribution $\eta := \sum V_i \delta_{X_i}$ defines a Dirichlet process with underlying parameter measure $\theta\nu$, where $\nu$ is the common law of the $X_i$'s. In this sense, $\mathrm{PD}(0, \theta)$ is the simplicial part of a Dirichlet process. By replacing $\mathrm{PD}(0, \theta)$ by $\mathrm{PD}(\alpha, \theta)$, a two-parameter generalization of a Dirichlet process was introduced in [49].

In the case of Dirichlet processes, there are many articles, including [5, 11, 23, 36, 46], in which exact forms (the distribution function, for example) of the law of the (random) mean $\sum X_i V_i$ of $\eta$ are obtained in terms of $\nu$. A key tool is an integral identity due to Cifarelli and Regazzini [5] which connects these two laws. As explained in [31, 50, 53], it is called the *Markov–Krein identity* because in the case $\theta = 1$, an integral transform is involved, analogous to the one studied by Markov and Krein in the context of moment problems; see [30] for background and various applications. Its extension was proven by Tsilevich [49] for the generalized Dirichlet process; see also [31] for a further extension and [50, 53] for a simple proof under some restriction on the support of $\nu$. This identity gives a one-to-one correspondence between $\nu$ and the law of $\sum X_i V_i$ with $(V_i)$ having $\mathrm{PD}(\alpha, \theta)$ distribution. However, such a correspondence is implicit and seems very subtle in general. One therefore needs some procedure of inversion in order to obtain explicit information. Recently, James, Lijoi, and Prünster [25] obtained some distributional results for the generalized Dirichlet process by means of a Perron–Stieltjes-type inversion formula. Since $\mathrm{PD}(0, \theta)$ and $\mathrm{PD}(\alpha, 0)$ are related to gamma processes and stable subordinators, respectively [45], the corresponding problem is naturally considered for a more general class of random distributions derived from subordinators, as posed in [53]. For results of this kind, we refer the reader to works by Regazzini, Lijoi and Prünster [47] and by Nieto-Barajas, Prünster and Walker [38], both of which make essential use of the Gurland inversion formula.

Our focus now will be on what is implied by Theorem 3.2 in the aforementioned context. It will be pointed out that the basic identity (3.8) involving the probability generating functional exhibits a mathematical structure underlying these kinds of identities. This should be compared with the proofs in [11, 31, 49], where the Markov–Krein identity is viewed as a relation between moment sequences through the Ewens(–Pitman) sampling formula [13, 42]. Let us introduce some notation used to describe the domain of a map defined via the Markov–Krein identity. Denoting by $\mathcal{P}$ the totality of Borel probability measures on $\mathbf{R}$, let

$$\mathcal{P}_0 = \left\{ \nu \in \mathcal{P} : \int_{\mathbf{R}} \nu(\mathrm{d}x) \log(1+|x|) < \infty \right\}$$

and set $\mathcal{P}_\alpha = \{\nu \in \mathcal{P} : \int_{\mathbf{R}} \nu(\mathrm{d}x)|x|^\alpha < \infty\}$ for $0 < \alpha < 1$. The Markov–Krein identity already mentioned is the following. Taking a $\nu \in \mathcal{P}_\alpha$, let $\{X_i\}_{i=1}^\infty$ be i.i.d. random variables with each law of $X_i$ being $\nu$ and let $(V_i)_{i=1}^\infty$ be $\mathrm{PD}(\alpha, \theta)$-distributed. If $\{X_i\}_{i=1}^\infty$ and $(V_i)_{i=1}^\infty$ are mutually independent, then $\sum |X_i| V_i < \infty$ a.s. (as was shown in [14] for $\alpha = 0$ and in Proposition 1 of [47] for $0 < \alpha < 1$) and the law, denoted $\mathcal{M}_{\alpha,\theta}\nu$, of $M := \sum X_i V_i$ is characterized by one of the following equalities:

(i) for $\alpha = 0$ and $\theta > 0$,

$$\int_{\mathbf{R}} \mathcal{M}_{0,\theta}\nu(\mathrm{d}x)(z-x)^{-\theta} = \exp\left(-\theta \int_{\mathbf{R}} \nu(\mathrm{d}x)\log(z-x)\right) \qquad (z \in \mathbf{C} \setminus \mathbf{R}); \tag{7.1}$$

(ii) for $\alpha \in (0,1)$ and $\theta > 0$,

$$\int_{\mathbf{R}} \mathcal{M}_{\alpha,\theta}\nu(\mathrm{d}x)(z-x)^{-\theta} = \left(\int_{\mathbf{R}} \nu(\mathrm{d}x)(z-x)^\alpha\right)^{-\theta/\alpha} \qquad (z \in \mathbf{C} \setminus \mathbf{R}); \tag{7.2}$$

(iii) for $\alpha \in (0,1)$ and $\theta \in (-\alpha, 0)$, $\mathcal{M}_{\alpha,\theta}\nu \in \mathcal{P}_{-\theta}$ and (7.2) holds true;

(iv) for $\alpha \in (0,1)$ and $\theta = 0$, $\mathcal{M}_{\alpha,0}\nu \in \mathcal{P}_0$ and

$$\int_{\mathbf{R}} \mathcal{M}_{\alpha,0}\nu(\mathrm{d}x)\log(z-x) = \alpha^{-1}\log\int_{\mathbf{R}} \nu(\mathrm{d}x)(z-x)^\alpha \qquad (z \in \mathbf{C} \setminus \mathbf{R}). \tag{7.3}$$

Thus the correspondence $\nu \mapsto \mathcal{M}_{\alpha,\theta}\nu$ defines a map from $\mathcal{P}_\alpha$ to $\mathcal{P}$, for which we also write $\mathcal{M}_{\alpha,\theta}$, by a slight abuse of notation. Since we are not aware of any reference in which integrability of the transformed measure $\mathcal{M}_{\alpha,\theta}\nu$ claimed in (iii) or (iv) has been shown, the proof shall be given below by applying our result. As will be seen later (Proposition 7.2(i)), such a property is also needed for further discussion. In the subsequent argument, we often use the following three equalities:

$$\int_0^\infty \mathrm{d}s\, s^{\theta-1}\mathrm{e}^{-s}\mathrm{e}^{-us} = \Gamma(\theta)(1+u)^{-\theta} \qquad (\theta > 0), \tag{7.4}$$

$$C_\alpha \int_0^\infty \mathrm{d}s\, s^{-\alpha-1}\mathrm{e}^{-s}(\mathrm{e}^{-us}-1) = 1-(1+u)^\alpha \qquad (0 < \alpha < 1), \tag{7.5}$$

$$\int_0^\infty \mathrm{d}s\, s^{-1}\mathrm{e}^{-s}(\mathrm{e}^{-us}-1) = -\log(1+u), \tag{7.6}$$

where $u \in \mathbf{C}$ is such that $\mathrm{Re}\, u > -1$.

**Lemma 7.1.** *Let $\alpha \in (0,1)$ and $\theta \in (-\alpha, 0]$. Then $\mathcal{M}_{\alpha,\theta}\nu \in \mathcal{P}_{-\theta}$ for any $\nu \in \mathcal{P}_\alpha$.*

**Proof.** Define $M^* := \sum |X_i| V_i$. By the assumed independence,

$$E[e^{-sM^*}] = E_{\alpha,\theta}\left[\prod_{i=1}^{\infty} \int_{\mathbf{R}} \nu(\mathrm{d}x) e^{-sV_i|x|}\right] \qquad (s \geq 0). \tag{7.7}$$

We will apply Theorem 3.2 by taking $g(s) = \int_{\mathbf{R}} \nu(\mathrm{d}x) \exp(-s|x|)$. For this purpose, with the help of Fubini's theorem, observe from (7.5) that for $0 < \alpha < 1$,

$$\frac{C_\alpha}{\lambda^\alpha} \int_0^\infty \mathrm{d}z\, \frac{e^{-\lambda z}}{z^{\alpha+1}} |g(z) - 1| = \int_{\mathbf{R}} \nu(\mathrm{d}x)(1 + \lambda^{-1}|x|)^\alpha - 1 < \infty, \tag{7.8}$$

which implies that $\lambda_\alpha^*(g) = 0$. Therefore, (3.8) holds for all $\lambda > 0$ and reads

$$\frac{\lambda^\theta}{\Gamma(\theta+1)} \int_0^\infty \mathrm{d}s\, s^{\theta-1} e^{-\lambda s} (E[e^{-sM^*}] - 1)$$

$$= \begin{cases} \dfrac{1}{\theta}\left(\displaystyle\int_{\mathbf{R}} \nu(\mathrm{d}x)(1 + \lambda^{-1}|x|)^\alpha\right)^{-\theta/\alpha} - \dfrac{1}{\theta}, & 0 < \alpha < 1,\ \theta \neq 0, \\[2ex] -\dfrac{1}{\alpha}\log\displaystyle\int_{\mathbf{R}} \nu(\mathrm{d}x)(1 + \lambda^{-1}|x|)^\alpha, & 0 < \alpha < 1,\ \theta = 0. \end{cases} \tag{7.9}$$

The desired integrability can be seen as follows. For $\theta = 0$, the left-hand side of (7.9) with $\lambda = 1$ becomes $-E[\log(1 + M^*)]$ by (7.6), while the right-hand side is finite for $\nu \in \mathcal{P}_\alpha$. This proves that $\mathcal{M}_{\alpha,0}\nu \in \mathcal{P}_0$ whenever $\nu \in \mathcal{P}_\alpha$. Similarly, for $\theta \in (-\alpha, 0)$, by using (7.5) with $-\theta$ in place of $\alpha$, we see that the left-hand side of (7.9) with $\lambda = 1$ equals $(E[(1 + M^*)^{-\theta}] - 1)/\theta$. Therefore, $\mathcal{M}_{\alpha,\theta}\nu \in \mathcal{P}_{-\theta}$ is implied by $\nu \in \mathcal{P}_\alpha$. $\qquad\square$

**Remarks.** (i) In view of (7.9), the reader will be able to see that the dichotomy result of Feigin–Tweedie type [14] (originally shown for PD$(0,\theta)$'s) holds true for PD$(\alpha,\theta)$ with $\theta > 0$. More precisely, if $\nu \in \mathcal{P} \setminus \mathcal{P}_\alpha$ and $M^*$ is as above, then $M^* = \infty$ a.s. Thus the maximal domain of $\mathcal{M}_{\alpha,\theta}$ with $\theta > 0$ is identified with $\mathcal{P}_\alpha$. In the Dirichlet process case, Cifarelli and Regazzini [6] gave a proof of this fact by using (7.1). Their cutoff argument also applies for $0 < \alpha < 1$, as seen in the following. Considering $M_n^* := \sum |X_i| \mathbf{1}_{\{|X_i| \leq n\}} V_i$ for each $n = 1, 2, \ldots$, which corresponds to the image measure $\nu_n$ of $\nu$ under $x \mapsto |x| \mathbf{1}_{[-n,n]}(x)$, we have (7.9) with $M_n^*$ and $\nu_n$ in place of $M^*$ and $\nu$, respectively. Under the assumption that $\nu \notin \mathcal{P}_\alpha$, letting $n \to \infty$ yields $E[(1 + M^*/\lambda)^{-\theta}] = 0$ and hence $M^* = \infty$ a.s., as long as $\theta > 0$.

(ii) The above proof makes it almost obvious that the Markov–Krein identity itself can be recovered from Theorem 3.2 with $g(s) = \psi_\nu(\pm s)$, where $\psi_\nu$ is the characteristic function of $\nu$. These choices are allowed for any $\nu \in \mathcal{P}_\alpha$ because for all $x \in \mathbf{R}$,

$$\int_0^\infty \mathrm{d}z\, \frac{e^{-z}}{z^{\alpha+1}} |e^{\sqrt{-1}xz} - 1| \leq \begin{cases} A_0\log(1 + |x|) + B_0\min\{|x|, 1\}, & \alpha = 0, \\ A_\alpha|x|^\alpha + B_\alpha\min\{|x|, 1\}, & 0 < \alpha < 1, \end{cases}$$

where $A_\alpha$ and $B_\alpha$ are constants depending only on $\alpha$. (For the proof of this inequality, first observe that the left-hand side is bounded above by $|x|\Gamma(1 - \alpha)$ and then show it for $x > 1$ by, for example, dividing the domain of integration into $(0, 1/x)$ and $[1/x, \infty)$.) In addition, (7.5) and (7.6) show that $\lambda_\alpha^*(\psi_\nu(\pm \cdot)) = 0$. Therefore, analogously to (7.7),

$$\psi_{\mathcal{M}_{\alpha,\theta}\nu}(\pm s) = E[\exp(\pm \sqrt{-1}sM)] = E_{\alpha,\theta}\left[\prod_{i=1}^\infty \psi_\nu(\pm sV_i)\right], \quad (7.10)$$

and we obtain, by (3.8),

$$\lambda^\theta \int_0^\infty ds\, s^{\theta-1} e^{-\lambda s} (\psi_{\mathcal{M}_{\alpha,\theta}\nu}(\pm s) - 1) = R_{\alpha,\theta}\left(\lambda^{-\alpha} \int_0^\infty dz\, \frac{e^{-\lambda z}}{z^{\alpha+1}}(\psi_\nu(\pm z) - 1)\right) \quad (7.11)$$

for all $\lambda > 0$. In the same way as before, (7.11) can be converted to (7.1)–(7.3) with $z = \pm\sqrt{-1}\lambda$ by making use of (7.4)–(7.6).

The rest of the paper is devoted to describing various consequences of the above calculations and the Markov–Krein identity.

**Proposition 7.2.** (i) *Let $0 < \beta < \alpha < 1$ and $\theta > -\beta$. Then, as a map defined on $\mathcal{P}_\alpha$,*

$$\mathcal{M}_{\beta,\theta} \circ \mathcal{M}_{\alpha,-\beta} = \mathcal{M}_{\alpha,\theta}. \quad (7.12)$$

(ii) *Let $0 \le \alpha < 1$ and $\theta > -\alpha$. Suppose that $\nu \in \mathcal{P}_\alpha$. Then, for a.e. $t \in \mathbf{R}$,*

$$\psi_{\mathcal{M}_{\alpha,\theta}\nu}(t) = 1 + \sum_{n=1}^\infty \frac{c_{n,\alpha,\theta}}{n!} \int_{\Delta_n} \prod_{i=1}^n \frac{\psi_\nu(tv_i) - 1}{v_i^{\alpha+1}}\left(1 - \sum_{j=1}^n v_j\right)^{\theta+\alpha n-1} dv_1 \cdots dv_n. \quad (7.13)$$

(iii) *Let $0 < \alpha < 1$ and $\theta > -\alpha$. Then $\mathcal{M}_{\alpha,\theta}$ defined on $\mathcal{P}_\alpha$ is injective.*

**Proof.** (i) First, note that by Lemma 7.1, the composition $\mathcal{M}_{\beta,\theta} \circ \mathcal{M}_{\alpha,-\beta}$ is well defined on $\mathcal{P}_\alpha$ under the conditions assumed. (7.12) with $\theta = 0$ becomes $\mathcal{M}_{\beta,0} \circ \mathcal{M}_{\alpha,-\beta} = \mathcal{M}_{\alpha,0}$, the proof of which requires careful handling of branches. Taking an arbitrary $\nu \in \mathcal{P}_\alpha$ and putting $\nu' = (\mathcal{M}_{\beta,0} \circ \mathcal{M}_{\alpha,-\beta})\nu$, we observe from (7.2) and (7.3) that

$$\int_{\mathbf{R}} \nu'(dx)\log(z-x) = \frac{1}{\beta}\log \int_{\mathbf{R}} \mathcal{M}_{\alpha,-\beta}\nu(dx)(z-x)^\beta$$

$$= \frac{1}{\beta}\log\left(\int_{\mathbf{R}} \nu(dx)(z-x)^\alpha\right)^{\beta/\alpha} = \frac{1}{\beta}\log\exp\left(\frac{\beta}{\alpha}\log\int_{\mathbf{R}} \nu(dx)(z-x)^\alpha\right).$$

It should here be noted that since $0 < \beta/\alpha < 1$,

$$\arg\exp\left(\sqrt{-1}\frac{\beta}{\alpha}\arg\int_{\mathbf{R}} \nu(dx)(z-x)^\alpha\right) = \frac{\beta}{\alpha}\arg\int_{\mathbf{R}} \nu(dx)(z-x)^\alpha.$$

Combining these equalities yields

$$\int_{\mathbf{R}} \nu'(\mathrm{d}x) \log(z-x) = \frac{1}{\beta} \cdot \frac{\beta}{\alpha} \log \int_{\mathbf{R}} \nu(\mathrm{d}x)(z-x)^{\alpha} = \int_{\mathbf{R}} \mathcal{M}_{\alpha,0}\nu(\mathrm{d}x) \log(z-x),$$

where the last equality follows from (7.3). Since $z \in \mathbf{C} \setminus \mathbf{R}$ is arbitrary, this implies $\nu' = \mathcal{M}_{\alpha,0}\nu$, as required. It only remains to deal with the case where $0 < \beta < \alpha < 1$ and $\theta \neq 0$ is such that $\theta > -\beta$. Calculations are quite similar to those in the previous case and so are left to the reader. The proof of (i) has already been done.

(ii) By virtue of (7.11), (7.13) is a direct consequence of Lemma 3.1(i).

(iii) First, suppose that $\mathcal{M}_{\alpha,0}\nu = \mathcal{M}_{\alpha,0}\nu'$ for some $\nu, \nu' \in \mathcal{P}_{\alpha}$. Then, by (7.3), $\int_{\mathbf{R}} \nu(\mathrm{d}x)(1+\sqrt{-1}\lambda^{-1}x)^{\alpha} = \int_{\mathbf{R}} \nu'(\mathrm{d}x)(1+\sqrt{-1}\lambda^{-1}x)^{\alpha}$ for all $\lambda \in \mathbf{R} \setminus \{0\}$. This, together with (7.5), implies that $\nu = \nu'$. Next, let $\theta \neq 0$ be such that $\theta > -\alpha$ and assume that $\mathcal{M}_{\alpha,\theta}\nu = \mathcal{M}_{\alpha,\theta}\nu'$ for some $\nu, \nu' \in \mathcal{P}_{\alpha}$. Taking care with branches, we see by (7.2) that for each $\lambda \in \mathbf{R} \setminus \{0\}$, there exists an integer $n(\lambda)$ such that

$$\log \int_{\mathbf{R}} \nu(\mathrm{d}x)(1+\sqrt{-1}\lambda^{-1}x)^{\alpha} - \log \int_{\mathbf{R}} \nu'(\mathrm{d}x)(1+\sqrt{-1}\lambda^{-1}x)^{\alpha} = 2\pi n(\lambda)\alpha\theta^{-1}. \quad (7.14)$$

However, since the left-hand side is continuous in $\lambda$ and tends to 0 as $|\lambda| \to \infty$, we have $n(\lambda) \equiv 0$. This proves $\nu = \nu'$ and therefore $\mathcal{M}_{\alpha,\theta}$ is injective.  $\square$

Let us make some comments on Proposition 7.2. The assertion (i) extends Theorem 2.1 of [25], where (7.12) with $\beta = 0$ is shown and applied for a sampling procedure. We can also understand (7.12) via (7.11) combined with $R_{\beta,\theta} \circ R_{\alpha,-\beta} = R_{\alpha,\theta}$. In general, a random variable $M$ having $(\mathcal{M}_{\alpha,\beta} \circ \mathcal{M}_{\gamma,\delta})\nu$ distribution (if such a law is well defined) is constructed by $M = \sum_{j} (\sum_{i} X_{ij} V_{ij}^{(\gamma,\delta)}) V_{j}^{(\alpha,\beta)}$, where $\{X_{ij}\}$ is a family of i.i.d. random variables with common law $\nu$, $(V_{j}^{(\alpha,\beta)})_{j=1}^{\infty}$ has PD$(\alpha,\beta)$ distribution and $(V_{i1}^{(\gamma,\delta)})_{i=1}^{\infty}$, $(V_{i2}^{(\gamma,\delta)})_{i=1}^{\infty}, \ldots$ are PD$(\gamma,\delta)$-distributed. Moreover, these random elements are required to be mutually independent. (ii) and (iii) are motivated by some results of Lijoi and Regazzini [36] concerned with $\mathcal{M}_{0,\theta}$. Applying to (7.11) a general complex inversion formula for Laplace transforms (the Bromwich integral), the authors provided an expression for the characteristic function of $\mathcal{M}_{0,\theta}\nu$. It is not clear how to verify equivalence between such a formula and (7.13) directly. They also proved a prototype of (iii), that is, that $\mathcal{M}_{0,\theta} : \mathcal{P}_0 \to \mathcal{P}$ is injective for each $\theta > 0$.

As observed [26] (Remark 4.2) in a much more general scheme called a *species sampling model*, Cauchy distributions are fixed points of every $\mathcal{M}_{\alpha,\theta}$. Alternatively, this also follows from (7.10). It should be emphasized that any property of PD$(\alpha,\theta)$ as a law on $\nabla_{\infty}$ is irrelevant to the validity of this fact, except the one which ensures that $\sum |X_i|V_i < \infty$ a.s. The study of the inverse problem would require a deeper understanding of $\mathcal{M}_{\alpha,\theta}$ itself or PD$(\alpha,\theta)$. In this respect, Yamato [55] gave a partial result concerning the case $\alpha = 0$. It was shown in [36] that for any $\nu \in \mathcal{P}_0$, $\mathcal{M}_{0,1}\nu = \nu$ occurs only when $\nu$ is Cauchy or degenerate. It seems reasonable to expect that such an equivalence would hold true for any $\mathcal{M}_{\alpha,\theta}$ and we could approach this problem by the use of results presented in Section

**3.** But here we shall show only a minor extension by mimicking the proof of Theorem 5 in [36], that is, by reducing to an ordinary differential equation via the Markov–Krein identity.

**Proposition 7.3.** *Let $0 < \alpha < 1$ and $\nu \in \mathcal{P}_\alpha$. Then $\mathcal{M}_{\alpha,1-\alpha}\nu = \nu$ if and only if $\nu$ is Cauchy or degenerate.*

**Proof.** Supposing that $\mathcal{M}_{\alpha,1-\alpha}\nu = \nu$, we only have to show that $\nu$ is Cauchy or degenerate. By combining this equality with (7.2), we easily see that $f_\pm(t) := \int_{\mathbf{R}} \nu(\mathrm{d}x)(t \pm \sqrt{-1}x)^\alpha$ solves the equation $\alpha^{-1} f_\pm{}'(t) = f_\pm(t)^{1-1/\alpha} (t > 0)$ and therefore $f_\pm(t)^{1/\alpha} = t + f_\pm(0)^{1/\alpha} (t \geq 0)$. Here, since $f_+(0)$ and $f_-(0)$ are complex conjugate with $|\arg f_\pm(0)| \leq \pi\alpha/2$, we have $|\arg(f_\pm(0)^{1/\alpha})| \leq \pi/2$ so that $f_\pm(0)^{1/\alpha} = \sigma^2 \pm \sqrt{-1}m$ for some $\sigma^2 \geq 0$ and $m \in \mathbf{R}$. Consequently, $f_\pm(t)^{1/\alpha} = t + \sigma^2 \pm \sqrt{-1}m (t \geq 0)$. Furthermore, by noting the continuity of $\arg f_\pm(t)$, we get $f_\pm(t) = (t + \sigma^2 \pm \sqrt{-1}m)^\alpha$ for every $t \geq 0$. With the help of (7.5), we conclude that $\nu$ is a Cauchy or degenerate distribution with $\psi_\nu(t) = \exp(\sqrt{-1}mt - \sigma^2|t|)$. $\qquad\square$

# Acknowledgements

This work was supported in part by the exchange program between the Japan Society for the Promotion of Science and Centre National de la Recherche Scientifique. It was initiated during a stay at Laboratoire de Probabilité et Modèles Aléatoires, Paris. The author would like to express his sincere gratitude to Professor Marc Yor for stimulating discussions and gratefully acknowledges the warm hospitality of the laboratory. He wishes to thank the anonymous referee for providing comments and help in improving the presentation of this paper, and is also grateful to S. Feng, who informed him of the recent article [16]. This paper is dedicated to Professor Itaru Mitoma on the occasion of his 60th birthday.

# References

[1] Arratia, R., Barbour, A.D. and Tavaré, S. (2003). *Logarithmic Combinatorial Structures: A Probabilistic Approach.* Zurich: European Mathematical Society. MR2032426

[2] Billingsley, P. (1995). *Probability and Measure*, 3rd ed. New York: Wiley. MR1324786

[3] Billingsley, P. (1999). *Convergence of Probability Measures*, 2nd ed. New York: Wiley. MR1700749

[4] Buchstab, A.A. (1937). An asymptotic estimation of a general number-theoretic function. *Mat. Sb.* **44** 1239–1246.

[5] Cifarelli, D.M. and Regazzini, E. (1990). Distribution functions of means of a Dirichlet process. *Ann. Statist.* **18** 429–442; correction. ibid. **22** (1994) 1633–1634. MR1041402

[6] Cifarelli, D.M. and Regazzini, E. (1996). Tail behaviour and finiteness of means of distributions chosen from a Dirichlet process. Quaderno IAMI 96.19. CNR-IAMI, Milano. Available at http://www.mi.imati.cnr.it/iami/abstracts/96-19.html.

[7] Daley, D.J. and Vere-Jones, D. (2003). *An Introduction to the Theory of Point Processes. Vol. I. Elementary Theory and Methods*, 2nd ed. New York: Springer. MR1950431

[8] Darling, D.A. (1952). The influence of the maximum term in the addition of independent random variables. *Trans. Amer. Math. Soc.* **73** 95–107. MR0048726

[9] de Bruijn, N.G. (1981). *Asymptotic Methods in Analysis*. New York: Dover Publications. (Corrected reprint of the third edition.) MR0671583

[10] Derrida, B. (1997). From random walks to spin glasses. *Phys. D* **107** 186–198. MR1491962

[11] Diaconis, P. and Kemperman, J. (1996). Some new tools for Dirichlet priors. In *Bayesian Statistics 5, (Alicante, 1994)* 97–106. New York: Oxford Univ. Press. MR1425401

[12] Dickman, K. (1930). On the frequency of numbers containing prime factors of a certain relative magnitude. *Ark. Mat. Astr. Fys.* **22** 1–14.

[13] Ewens, W.J. (1972). The sampling theory of selectively neutral alleles. *Theoret. Population Biology* **3** 87–112; erratum. ibid. **3** (1972) 240, 376. MR0325177

[14] Feigin, P.D. and Tweedie, R.L. (1989). Linear functionals and Markov chains associated with Dirichlet processes. *Math. Proc. Cambridge Philos. Soc.* **105** 579–585. MR0985694

[15] Feller, W. (1968). *An Introduction to Probability Theory and Its Applications. Vol. I*, 3rd ed. New York: Wiley. MR0228020

[16] Feng, S. (2007). Large deviations for Dirichlet processes and Poisson–Dirichlet distribution with two parameters. *Electron. J. Probab.* **12** 787–807. MR2318410

[17] Ferguson, T.S. (1973). A Bayesian analysis of some nonparametric problems. *Ann. Statist.* **1** 209–230. MR0350949

[18] Griffiths, R.C. (1979). On the distribution of allele frequencies in a diffusion model. *Theoret. Population Biol.* **15** 140–158. MR0528914

[19] Griffiths, R.C. (1988). On the distribution of points in a Poisson Dirichlet process. *J. Appl. Probab.* **25** 336–345. MR0938197

[20] Grote, M.N. and Speed, T.P. (2002). Approximate Ewens formulae for symmetric overdominance selection. *Ann. Appl. Probab.* **12** 637–663. MR1910643

[21] Hensley, D. (1984). The sum of $\alpha^{\Omega(n)}$ over integers $n \leq x$ with all prime factors between $\alpha$ and $y$. *J. Number Theory* **18** 206–212. MR0741951

[22] Ignatov, T. (1982). On a constant arising in the asymptotic theory of symmetric groups, and Poisson–Dirichlet measures. *Theory Probab. Appl.* **27** 136–147. MR0645134

[23] James, L.F. (2005). Functionals of Dirichlet processes, the Cifarelli–Regazzini identity and beta–gamma processes. *Ann. Statist.* **33** 647–660. MR2163155

[24] James, L.F. (2008). Large sample asymptotics for the two-parameter Poisson–Dirichlet process. *IMS Collections* **3** 187–200. MR2459225

[25] James, L.F., Lijoi, A. and Prünster, I. (2008). Distributions of linear functionals of two parameter Poisson–Dirichlet random measures. *Ann. Appl. Probab.* **18** 521–551. MR2398765

[26] James, L.F., Lijoi, A. and Prünster, I. (2009). On the posterior distribution of means of nonparametric prior processes. *Bernoulli*. To appear.

[27] Johansson, K. (2005). Random matrices and determinantal processes. Lecture notes from the Les Houches summer school on mathematical statistical physics. Available at arXiv:math-ph/0510038.

[28] Joyce, P., Krone, S.M. and Kurtz, T.G. (2002). Gaussian limits associated with the Poisson–Dirichlet distribution and the Ewens sampling formula. *Ann. Appl. Probab.* **12** 101–124. MR1890058

[29] Joyce, P., Krone, S.M. and Kurtz, T.G. (2003). When can one detect overdominant selection in the infinite-alleles model? *Ann. Appl. Probab.* **13** 181–212. MR1951997

[30] Kerov, S. (1998). Interlacing measures. In *Kirillov's Seminar on Representation Theory. Amer. Math. Soc. Transl. Ser. 2* **181** 35–83. Providence, RI: Amer. Math. Soc. MR1618739

[31] Kerov, S.V. and Tsilevich, N.V. (2001). The Markov–Krein correspondence in several dimensions. *Zap. Nauchn. Sem. POMI* **283** 98–122; translation. *J. Math. Sci. (New York)* **121** (2004) 2345–2359. MR1879065

[32] Kingman, J.F.C. (1975). Random discrete distribution. *J. Roy. Statist. Soc. Ser. B* **37** 1–22. MR0368264

[33] Kingman, J.F.C. (1993). *Poisson Processes.* New York: Oxford Univ. Press. MR1207584

[34] Kolchin, V.F. (1986). *Random Mappings.* New York: Optimization Software. MR0865130

[35] Lamperti, J. (1961). A contribution to renewal theory. *Proc. Amer. Math. Soc.* **12** 724–731. MR0125663

[36] Lijoi, A. and Regazzini, E. (2004). Means of a Dirichlet process and multiple hypergeometric functions. *Ann. Probab.* **32** 1469–1495. MR2060305

[37] McCloskey, J.W. (1965). A model for the distribution of individuals by species in an environment. Ph.D. thesis, Michigan State University.

[38] Nieto-Barajas, L.E., Prünster, I. and Walker, S.G. (2004). Normalized random measures driven by increasing additive processes. *Ann. Statist.* **32** 2343–2360. MR2153987

[39] Olshanski, G. (1998). Point processes and the infinite symmetric group. Part I: The general formalism and the density function. Available at arXiv:math/9804086v1. MR1671191

[40] Penrose, M.D. and Wade, A.R. (2004). Random minimal directed spanning trees and Dickman-type distributions. *Adv. in Appl. Probab.* **36** 691–714. MR2079909

[41] Perman, M., Pitman, J. and Yor, M. (1992). Size-biased sampling of Poisson point processes and excursions. *Probab. Theory Related Fields* **92** 21–39. MR1156448

[42] Pitman, J. (1995). Exchangeable and partially exchangeable random partitions. *Probab. Theory Related Fields* **102** 145–158. MR1337249

[43] Pitman, J. (1996). Random discrete distributions invariant under size-biased permutation. *Adv. in Appl. Probab.* **28** 525–539. MR1387889

[44] Pitman, J. (2006). *Combinatorial Stochastic Processes. Lecture Notes in Mathematics* **1875**. Berlin: Springer. MR2245368

[45] Pitman, J. and Yor, M. (1997). The two-parameter Poisson–Dirichlet distribution derived from a stable subordinator. *Ann. Probab.* **25** 855–900. MR1434129

[46] Regazzini, E., Guglielmi, A. and Di Nunno, G. (2002). Theory and numerical analysis for exact distributions of functionals of a Dirichlet process. *Ann. Statist.* **30** 1376–1411. MR1936323

[47] Regazzini, E., Lijoi, A. and Prünster, I. (2003). Distributional results for means of normalized random measures with independent increments. *Ann. Statist.* **31** 560–585. MR1983542

[48] Tennenbaum, G. (1995). *Introduction to Analytic and Probabilistic Number Theory.* London: Cambridge Univ. Press. MR1342300

[49] Tsilevich, N. (1997). Distribution of mean values for some random measures. *Zap. Nauchn. Sem. POMI* **240** 268–279; translation *J. Math. Sci. (New York)* **96** (1999) 3616–3623. MR1691650

[50] Tsilevich, N., Vershik, A. and Yor, M. (2000). Distinguished properties of the gamma process, and related topics. *Prépublication du Laboratoire de Probabilités et Modèles Aléatoires* No. 575. Available at http://xxx.lanl.gov/ps/math.PR/0005287.

[51] Tsilevich, N., Vershik, A. and Yor, M. (2001). An infinite-dimensional analogue of the Lebesgue measure and distinguished properties of the gamma process. *J. Funct. Anal.* **185** 274–296. MR1853759

[52] Vershik, A. and Yor, M. (1995). Multiplicativité du processus gamma et étude asymptotique des lois stables d'indice $\alpha$, lorsque $\alpha$ tend vers 0. Technical Report 289. Laboratoire de Probabilités, Univ. Paris VI.

[53] Vershik, A., Yor, M. and Tsilevich, N. (2001). On the Markov–Krein identity and quasi-invariance of the gamma process. *Zap. Nauchn. Sem. POMI* **283** 21–36; translation. *J. Math. Sci. (New York)* **121** (2004) 2303–2310. MR1879060

[54] Watterson, G.A. (1976). The stationary distribution of the infinitely-many neutral alleles diffusion model. *J. Appl. Probab.* **13** 639–651; correction. ibid. **14** (1977) 897. MR0504014

[55] Yamato, H. (1984). Characteristic functions of means of distributions chosen from a Dirichlet process. *Ann. Probab.* **12** 262–267. MR0723745